\documentclass[12pt]{amsart}
\usepackage{amscd,amssymb}

\usepackage[cmtip,matrix,arrow]{xy}
\numberwithin{equation}{section}

\newtheorem {Theorem} 			{Theorem}

\newtheorem {Lemma}[equation]     	{Lemma}
\newtheorem {Proposition}[equation]	{Proposition}

\newtheorem {corollary}[equation]	{Corollary}
\newtheorem {refproposition}		{Proposition}
\theoremstyle{definition}
\newtheorem {Definition}[equation]{Definition}
\newtheorem {Remark}[equation]		{Remark}
\newtheorem {Example}[equation]		{Example}

\newcommand{\pr} {\smallskip\noindent{\bf Proof\,\,}}

\newcommand     {\comment}[1]   {}
\newcommand     {\mute}[2] {}
\newcommand     {\printname}[1] {}

\newcommand{\labell}[1] {\label{#1}\printname{#1}}

\def	\hol	{{\operatorname{hol}}}

\def	\Stab	{\operatorname{Stab}}

\def	\inv	{^{-1}}
\def	\to	{\longrightarrow}

\def	\mod	{{/ \! /}}

\def	\orho	{{\overline {\rho}}}

\def	\R	{{\mathbb R}}

\def	\cO	{{\mathcal O}}

\def	\ft	{{\mathfrak t}}
\def	\fk	{{\mathfrak k}}
\def	\fh	{{\mathfrak h}}
\def	\fg	{{\mathfrak g}}

\def	\Ad	{{\operatorname{Ad}}}

\def	\stab	{{\operatorname{stab }}}

\def	\hX 	{{ \hat X }}
\def	\hmu 	{{ \hat \mu }}

\def\ogh{\widehat{\Omega G}}
\def\ogs{{\Omega G}^*}
\def\pgh{\widehat{P_e G}}

\def\pgeh{\widehat{P_e G}}
\def\pges{{P_e G}^*}

\def	\pr	{\operatorname{pr}}

\begin{document}

%%%%%%%%%%%%%%%%%%%%%%%%%%%%%%%%%%%%%%%%%%%%%%%%%%%%%%%%%%%%%%%%%%%%%%%%%%%%%%%
%%%%%%%%%%%%%%%%%%%%%%%%%%%%%%%%%%%%%%%%%%%%%%%%%%%%%%%%%%%%%%%%%%%%%%%%%%%%%%%

\title[Surjectivity for Hamiltonian Loop Group Spaces]
{Surjectivity for Hamiltonian Loop Group Spaces}

\author{Raoul Bott}
\author{Susan Tolman}
\author{Jonathan Weitsman}
\thanks{R. Bott was supported in part by NSF grant
DMS 01/04196.  S. Tolman was partially supported by an NSF grant and
by an Alfred P. Sloan Foundation fellowship.
J. Weitsman was supported in part by NSF grant DMS 99/71914, by NSF
Young Investigator grant DMS 94/57821, and
by an Alfred P. Sloan Foundation Fellowship.}
\address{Department of Mathematics, Harvard University, Cambridge,
MA  02138}
\email{bott@math.harvard.edu}
\address{Department of Mathematics,  University of Illinois at
Urbana-Champaign,
Urbana, IL 61801}
\email{stolman@math.uiuc.edu}
\address{Department of Mathematics, University of California, Santa Cruz,
CA 95064}
\email{weitsman@cats.ucsc.edu}
\thanks{\today}

\begin{abstract}
Let $G$ be a compact Lie group, and let $LG$ denote the corresponding
loop group.
Let $(X,\omega)$ be a weakly symplectic Banach manifold. 
Consider a Hamiltonian action of $LG$ on $(X,\omega)$, and assume that
the moment map $\mu: X \to L\fg^*$ is proper.
We consider the function $|\mu|^2: X \to \R$, and use a version
of Morse theory to show that the
inclusion map $j:\mu^{-1}(0)\to X$ induces
a surjection $j^*:H_G^*(X) \to H_G^*(\mu^{-1}(0))$, in analogy
with Kirwan's surjectivity theorem in the finite-dimensional case.
We also prove a version of this surjectivity 
theorem for quasi-Hamiltonian $G$-spaces.

\end{abstract}

\maketitle

%%%%%%%%%%%%%%%%%%%%%%%%%%%%%%%%%%%%%%%%%%%%%%%%%%%%%%%%%%%%%%%%%%%%%%%%%%%%%%
%%%%%%%%%%%%%%%%%%%%%%%%%%%%%%%%%%%%%%%%%%%%%%%%%%%%%%%%%%%%%%%%%%%%%%%%%%%%%%
\section{Introduction}

\comment{Integrate examples more effectively}

Let $G$ be a compact Lie group. 
A {\bf Hamiltonian $G$ space} is a triple $(M,\omega,\mu)$,
where $(M,\omega)$ is a symplectic manifold,
and $\mu: M \to \fg^*$ is a moment map for a Hamiltonian
action of $G$ on $M$.

If $0$ is a regular value of $\mu$, then
the {\bf symplectic quotient}
$M \mod G = \mu^{-1}(0)/G$  is a symplectic orbifold.
In this case, $G$ acts locally freely on $\mu^{-1}(0)$, so the
equivariant cohomology ring $H_G^*(\mu^{-1}(0))$ coincides
with the cohomology ring $H^*(M \mod G)$.\footnote{In this
paper, we consider only cohomology groups with rational coefficients.}
Therefore, the restriction
map from $H_G^*(M)$ to $H_G^*(\mu\inv(0))$ induces a
map $\kappa \colon H_G^*(M) \to H^*(M \mod G)$, which we call 
the {\bf Kirwan map}.

A version of Morse theory due to Kirwan \cite{kirwan} is
an important tool in the study of the topology of Hamiltonian
$G$-spaces and their symplectic quotients.  Kirwan shows that
for a compact Hamiltonian $G$-space $(M,\omega,\mu)$,
the function $|\mu|^2 \colon M \to \R$ 
is an equivariantly perfect Morse function on $M$.   This fact
has the following striking consequence.

\begin{Theorem} [Surjectivity for compact group actions (Kirwan)]
Let $G$ be a compact Lie group, and let $(M,\omega,\mu)$ be
a compact Hamiltonian $G$-space.
Then the restriction map
$H_G^*(M)\to H_G^*(\mu^{-1}(0))$ is surjective.
In particular, if $0$ is a regular value of $\mu$, the Kirwan map is
surjective.
\end{Theorem}

The main purpose of this paper is to generalize Theorem 1 to
the case of Hamiltonian loop group actions on symplectic
Banach manifolds with proper moment map.

Let $G$ be a compact Lie group, and let $\langle\cdot,\cdot\rangle$
be an invariant inner product on the Lie algebra $\fg$.
Let $L_s G$ be the space of maps from $S^1$ to
$G$ of Sobolev class $s> 1/2$; $L_s G$ is a Banach Lie group with
the group structure given
by pointwise multiplication.  The Lie algebra $L_s {\mathfrak g}$
is given by the space $\Omega_{s}^0(S^1) \otimes {\mathfrak g}$
of maps from $S^1$ to ${\mathfrak g}$ of Sobolev class $s$.
We define $L_s{\mathfrak g}^*$ to be the space
$\Omega_{s-1}^1(S^1) \otimes {\mathfrak g}$
of ${\mathfrak g}$-valued one forms on $S^1$ of Sobolev class $s-1$.
Integration gives a natural non-degenerate pairing $(\cdot, \cdot)$
of $L_s\fg$ with $L_s \fg^*$.
Using this inner product, the space
$L_s{\mathfrak g}^*$ can be identified with the space of
connections on the trivialized principal $G$-bundle
$G\times S^1\to S^1$.  This identification induces
an action of the group $L_s G$
on $L_s{\mathfrak g}^*$ given by\footnote
{
This action is not the coadjoint action of $L_sG$, 
but arises instead
from the coadjoint action of a central extension
of $L_s G$.  As a consequence it is an affine action
rather than a linear action on the vector space
$L_s{\mathfrak g}^*.$
}
\begin{equation}\labell{action}
\lambda \cdot \gamma  = \lambda^{-1} d \lambda + \lambda^{-1}\gamma \lambda
\end{equation}
where $\gamma \in L_s{\mathfrak g}^*$ and $\lambda \in L_s G$.

A Banach manifold $X$ is {\bf weakly symplectic}
if is is equipped with a closed two-form $\omega \in \Omega^2(X)$
such that the induced map  $\omega^\flat  \colon T_pX \to T_pX^*$
is injective.
An action of the group $L_sG$ on  $(X,\omega)$
is {\bf Hamiltonian} if there exists an $L_sG$-equivariant {\bf moment map}
$\mu  \colon X_s \to  L_s {\mathfrak g}^*$
so that $\langle d\mu, \xi\rangle = i_{\xi^X} \omega$ for all
$\xi \in L_s{\mathfrak g}$.
We call the triple
$(X,\omega, \mu)$ a {\bf Hamiltonian $L_s G$ space}\footnote
{
Let a Lie group $G$ act on a symplectic manifold $(M,\omega)$,
and let $\mu  \colon M \to  \fg*$
satisfy $\langle d\mu, \xi\rangle = i_{\xi^X} \omega$ for all
$\xi \in \fg$.  
If $G$ is compact, then one can choose 
$\mu$ to be equivariant
with respect to the coadjoint action of $G$ on $\fg^*$.
In general, however, this is not possible, and it
is natural  to consider  other actions of $G$ on $\fg^*$,
such as the action arising from a central extension in (\ref{action}).
}.

Now consider  a Hamiltonian $L_s G$ space $(X,\omega,\mu)$
with proper moment map. In this case,
if $0$ is a regular value of $\mu$,
the {\bf symplectic quotient}
$X \mod L_sG = \mu^{-1}(0)/G$ is a finite dimensional symplectic orbifold.
Since $G$ acts locally freely on $\mu^{-1}(0)$, the
equivariant cohomology ring $H_G^*(\mu^{-1}(0))$ coincides
with the cohomology ring $H^*(X \mod L_sG)$.
Therefore, the restriction map $H_G^*(X) \to H_G^*(\mu\inv(0))$
induces a  map $\kappa \colon H_G^*(X) \to H^*(X \mod L_sG)$
which is the  infinite-dimensional analog of the Kirwan map.

In this paper we consider the
the function $|\mu|^2 \colon  X\to \R$, and show that
it can be treated as an equivariantly
perfect Morse function.  We thus prove the following analog of
Kirwan's surjectivity theorem:

\def\dlg{L_s{\mathfrak g}^*}
\begin{Theorem}[Surjectivity for Hamiltonian loop group actions]
Let $G$ be a compact Lie group and choose 
an invariant inner product on $\fg$. Let $L_s G$ be the 
corresponding loop group, where $s > 1/2$.
Let $(X,\omega,\mu)$ be a Hamiltonian $L_s G$-space with
proper moment map.
Then the restriction $H_G^*(X) \to H_G^*(\mu^{-1}(0))$ is surjective.
In particular, if $0$ is a regular value of $\mu$, then
the Kirwan map $\kappa \colon H_G^*(X) \to H^*(X \mod L_sG)$ is surjective.
\end{Theorem}

\subsection{Proof of Theorem 1}

We begin by reviewing Kirwan's proof of surjectivity
in the finite-dimensional case.  Let $G$ be a compact
Lie group, and let $(M,\omega, \mu)$ be a Hamiltonian $G$-space
with proper moment map.
Consider the function $f = |\mu|^2 \colon  M \to \R$.
Since $\mu\inv(0) = f\inv(0)$, it 
is enough to prove that the restriction 
$H_G^*(M) \to H_G^*(f\inv(0))$ is
surjective.

The first step in Kirwan's argument is to show that
$f$ is {\bf Morse in the sense of Kirwan}
(see Definition \ref{morsekirwan}).
As a result,
for every component $C$ of the critical set of $f$
there exists a vector bundle $E_C^-$ over $C$, called the
{\bf negative normal bundle at $C$} (see Definition \ref{negnor}).
For sufficiently small $\epsilon > 0$, let
$M_\pm=f^{-1}(-\infty, f(C) \pm\epsilon)$.
Consider the long exact sequence of the pair
$(M_+, M_-)$ in equivariant cohomology.
By Proposition \ref{morsekirwanlemma}, we obtain a
commuting diagram

\begin{equation}\labell{exactseqintro}
\vcenter{\xymatrix{
\cdots
\ar[r] &
H_G^*(M_+,M_-)
\ar[r]\ar[d]^\simeq &
H_G^*(M_+)
\ar[r] \ar[d] &
H_G^*(M_-)\ar[r] &  \cdots
\\
&
H_G^{*-\lambda_C} (C)
\ar[r]^{\ \ \cup \, e_C}  &
H_G^*(C) &
\\
}}
\end{equation}
\noindent where $\lambda_C$ denotes the index of $C$ and
$e_C \in H_G^*(C)$ denotes the
equivariant Euler class of the negative normal bundle at $C$.

Kirwan now applies the completion principle of Atiyah and Bott \cite{ab}
to show that $f$ is an  equivariantly perfect Morse function.
The key idea is to show that the Euler classes
$e_C \in H_G^*(C)$ are not zero divisors.  By the diagram above,
this implies that the long exact sequence in equivariant relative cohomology
splits into short exact sequences

\begin{equation}\labell{shortexactseqintro}
0\to H_G^*(M_+,M_-) \to H^*_G(M_+) \to H^*_G(M_-) \to 0.
\end{equation}

\noindent
Since $f$ is non-negative, Theorem 1 follows by induction on the
critical levels.

To show that $e_C$ is not a zero divisor, Kirwan shows
that for each component $C$ of the critical set
there exists a subtorus $T \subset G$ and a $Z(T)$-invariant
subset $B \subset C^T$
$$G \times_{Z(T)} B \to C$$
is an equivariant
homeomorphism.\footnote{In this paper, given   
given a subgroup $K \subset G$, $Z(K)$ denotes the centralizer of $K$ in $G$.
If $G$ acts on a space $X$, then $X^K$ denotes the set of points
in $X$ fixed by $K$.}
Hence, we have a natural isomorphism
$H_G^*(C) \simeq H_{Z(T)}^*(B)$; this isomorphism takes the $G$-equivariant
Euler class of $E_C^-$ to the
$Z(T)$-equivariant Euler class of  $E_C^-|_{B}$.
Moreover, $(E_C^-)^T$ is a subset of the zero section of $E_C^-$.
The fact that $e_C$ is not a zero divisor now follows from the Lemma below.

\begin{Lemma}[Atiyah-Bott]\labell{atiyahbottlemmaintro}
Let $V$ be a complex vector bundle over a connected
space $Y$, and let a compact Lie group $K$ act on $V$.
If there exists a subtorus $T \subset K$
so that the fixed point set $V^T$ is the zero section of $V$, then
the Euler class $e(V) \in H_K^*(Y)$ is not a zero divisor.
\end{Lemma}

\subsection{Proof of Theorem 2.}

The basic idea of this proof is straightforward:
Morally, we would like to consider the function $f= |\mu|^2$ as
a Morse function on $X$, and follow the proof of Kirwan's surjectivity theorem  for
compact Hamiltonian $G$-spaces which we recalled in Section
1.2 above.
The main technical hurdle we must face is that Kirwan's extension
of Morse theory has only been developed in the finite dimensional case.
However, in the case where $X=L_sG/G$ is the Hamiltonian $L_sG$-space
given in Example \ref{Ex 1} below,
the function $f$ is precisely the classical energy
function of Morse and Bott \cite{Bott}.
In the spirit of \cite{Bott}, we will replace $X$ in the general case by
a sequence of finite dimensional approximations.

\def\hf{\hat{f}}
First, in Section 2,
we construct an {\bf infinite approximating space}
$\hat{X}$, along with an {\bf energy function}
$\hf \colon \hat{X} \to \R$,  and prove:

\begin{refproposition}[Proposition \ref{approxi}]
Let $G$ be a compact Lie group, and $L_s G$ be the corresponding
loop group, where $s > 1/2$.
Let $(X,\omega,\mu)$ be a Hamiltonian $L_s G$ space with
proper moment map.
Let $\hat{X}$ be the associated infinite approximating
space, and let $\hat{f} \colon \hat{X} \to \R$ be the energy function.

The restriction map
$H^*_G(X) \to H^*_G(f\inv(0))$ is surjective
if and only if the  restriction map
$ H^*_G(\hX) \to H^*_G(\hf\inv(0))$ is surjective.
\end{refproposition}

Next, in Section 3,
we construct a series of {\bf finite approximating spaces}
$Y_n$, each of which has an energy function $f_n \colon Y_n \to \R$.
We then prove:

\begin{refproposition}[Proposition \ref{finite}]
Let $(X,\omega,\mu)$ be a Hamiltonian $L_s G$ space with
proper moment map.
Let $\hX$ be the associated infinite approximating
space and let $\hf \colon \hX \to \R$ be the energy function.
Let $Y_n$ be the finite approximating manifolds,
and let $f_n \colon Y_n \to \R$ be the energy functions.

If the restriction maps $H^*_G(Y_n)\to H^*_G(f_n \inv(0))$ are all
surjections, then the restriction map  $H_G^*(\hat{X}) \to H_G^*(\hf
\inv(0))$
is a surjection.
\end{refproposition}

Since $f\inv(0) = \mu\inv(0)$, in order to prove
Theorem 2, it is enough to show that the restriction
$H_G^*(Y_n) \to H^*_G( f_n \inv(0))$ is surjective.
Unfortunately
the manifolds $Y_n$ are not symplectic, so that we cannot
directly apply Kirwan's results.
Instead, in Sections 4 through 8 of this
paper we work through the
local calculations needed to prove the two
results below.

\begin{refproposition}[Proposition \ref{fnisMK}]
The functions $f_n \colon Y_n \to \R$ are Morse in the sense of Kirwan.
\end{refproposition}

\begin{refproposition}[Proposition \ref{fnfixed}]
Let $C$ be a component of the critical set of $f_n$
and let $E_C^-$ be the negative normal bundle at $C$.
Then there exists a subtorus $T \subset G$ and a $Z(T)$ invariant
subset $B \subset C^T$
so that the natural map $G \times_{Z(T)} B \to C$ is an
equivariant homeomorphism.
Moreover,
$(E_C^-)^T$ is a subset of the zero section of $E_C^-$.
\end{refproposition}

We can now follow Kirwan's proof.
Proposition \ref{fnisMK} shows that for each
component of the critical set of $f_n$ we have a commuting diagram
analogous to
(\ref{exactseqintro}).
As in Kirwan's argument, Proposition \ref{fnfixed} now shows that
the equivariant Euler class of the negative normal
bundle at each component of the critical set is not a zero divisor.
Hence the long exact sequence in relative equivariant
cohomology breaks up into short exact sequences
as in (\ref{shortexactseqintro}).
Since each $f_n$ is non-negative, it follows by induction that
the restriction maps $H_G^*(Y_n)\to H_G^*(f_n \inv(0))$ are
surjections.

\subsection{Examples}
\def\dlgnos {L{\mathfrak g}^*}
The following are examples of Hamiltonian $LG$
spaces \cite{AMM}.\footnote{We suppress
the Sobolev index $s$ in $L_sG$ for convenience; the examples
in this section can be considered for any $s > 1/2$.}

\begin{Example}\label{Ex 1}[The minimal coadjoint orbit]
Consider $X=LG \cdot 0$, the $LG$-orbit of the point $0 \in \dlgnos$,
equipped with the Kirillov-Kostant symplectic structure
\cite{ps}.
The inclusion $X \to L\fg^*$ is a proper moment map for the
natural $LG$ action.
As $\Stab(0)=G$, this space may be identified
with the quotient $LG/G$.
In terms of the identification of $X$ with $LG/G$, the
symplectic form arises from the $G$-invariant two-form
$\tilde{\omega} \in \Omega^2(LG)$ given by
$$\omega|_{\lambda}(\xi,\eta) = {\frac 1 {2\pi}} \int_0^{2 \pi}
\langle \xi(t),\eta'(t)\rangle dt,$$
where $\xi,\eta \in TLG|_{\lambda}$ (see \cite{ps}, p. 147.).
The moment map for the $LG$ action is the map
${\mu} \colon LG/G\to \dlgnos$
arising from the $G$-invariant map
$\tilde{\mu} \colon LG\to \dlgnos$ defined by
$\tilde{\mu}(\lambda)= \lambda^{-1} d \lambda$.
The manifold $LG/G$ is diffeomorphic to
the subgroup $\Omega G = \{\lambda \in L G  \mid \lambda(0) = e\}$
of based loops in $L G.$ 
In terms of this identification, the square of the moment
map $E=|\mu|^2$ is the energy functional
$E(\lambda) =
{\frac 1 {2\pi}} \int_0^{2 \pi} |\lambda^{-1}{\frac {d\lambda} {dt}}|^2
dt$, studied by Morse and one of the present authors
(R.B.).  The perfection of $E$ as a Morse function plays
an important role in the proof of the Periodicity Theorem \cite{Bott}.
\end{Example}

\begin{Example}\label{Ex 2}[The generic coadjoint orbit]
Consider next a generic point $\xi \in \dlgnos$
and the corresponding coadjoint orbit $X=LG \cdot \xi$.  We may as well take
$\xi \in {\mathfrak t}^*$.  Then $\Stab(\xi)= T$, so that this space
may be identified with $LG/T$.  It it therefore diffeomorphic
to the space of paths in $G$ beginning at $e$ and
ending at a point of the conjugacy class of
$\exp \xi$.  Again, this is a Hamiltonian $LG$ space
with the Kirillov-Kostant symplectic form and with moment map
given by the inclusion $X\to L{\mathfrak g}^*$.
\end{Example}

\begin{Example}\label{Ex 3}[Spaces of connections on 2-manifolds]

Let $\Sigma$ be a compact, connected 2-manifold of genus $h$ with boundary
$\partial \Sigma =S^1$.  Consider the space
${\mathcal A}(\Sigma)= \Omega^1(\Sigma, {\mathfrak g})$
of connections on the trivialized principal $G$-bundle
$G\times \Sigma \to \Sigma$.  The space ${\mathcal A}(\Sigma)$
is a symplectic manifold, equipped with a Hamiltonian action
of the gauge group ${\mathcal G}(\Sigma)= Map(\Sigma,G)$.  The
moment map
$J \colon {\mathcal A}(\Sigma)\to Lie({\mathcal G}(\Sigma))^* =
\Omega^2(\Sigma,{\mathfrak g})$ is given
by $J(A) = F_A$, where $F_A$ is the curvature of the
connection $A$.  The group ${\mathcal G}(\Sigma)$ has a normal subgroup
${\mathcal G}(\Sigma,\partial \Sigma)$ defined by
${\mathcal G}(\Sigma,\partial \Sigma)= \{\gamma \in {\mathcal G}(\Sigma)
\mid
\gamma|_{\partial \Sigma}= e\}$,
and ${\mathcal G}(\Sigma)/{\mathcal G}(\Sigma,\partial \Sigma)=LG$.  Consider
the reduced space
${\mathcal A}(\Sigma) \mod {\mathcal G}(\Sigma,\partial \Sigma)$.  This is a
Hamiltonian $LG$-space with
proper moment map $\mu$.  The reduced space $X\mod LG$ is
the moduli space of flat connections on the trivial principal $G$-bundle
$G\times \Sigma \to \Sigma$.
\end{Example}

\mute
{The Yang-Mills functional $f=|F_A|^2$ on the space
of connections on a 2-manifold $\Sigma$ was studied
by Atiyah and Bott \cite{ab}.  The critical point
structure of the Yang-Mills functional is shown
there to give the correct Poincare polynomial for
the moduli space of flat connections.  Since the
group ${\mathcal G}(\Sigma, \partial\Sigma)$ acts
freely on ${\mathcal A}$, one would expect the
critical point structure of $|\mu|^2$ on $X$ to
be precisely the same as that of the Yang-Mills
functional, so that Theorem 2 should give a
new, infinite-dimensional proof of the results
of \cite{ab}.}

{\bf Acknowledgements:}  We would like to thank 
Anton Alexeev, Victor Guillemin, Eckhard Meinrenken,
Shlomo Sternberg and Chris Woodward (see \cite{W}) for helpful comments and
discussions.

\section{Surjectivity for quasi-Hamiltonian G-spaces}\labell{qham}

In this section we will review Alexe'ev, Malkin, and Meinrenken's
\cite{AMM}
definition of quasi-Hamiltonian $G$-spaces, together
with their proof that these spaces  are in one-to-one correspondence 
with Hamiltonian $L_s G$-spaces with proper moment map.
We then restate our surjectivity
theorem in terms of quasi-Hamiltonian $G$-spaces.

We begin by motivating the definition of quasi-Hamiltonian $G$-spaces.
Let $(X,\omega,\mu)$ be a Hamiltonian $L_s G$-space with proper
moment map.
The normal subgroup $\Omega_s G \subset L_sG$ given
by the based maps in $L_s G$
$$\Omega_s G = \{\lambda \in L_s G  \mid  \lambda(0) = e\}$$
acts freely on $L_s {\mathfrak g}^*$.
Since $\mu$ is $L_s G$-equivariant,  $\Omega_s G$ acts freely on $X$;
thus $X$ is a the
total space of a principal $\Omega_s G$ bundle $\pi \colon X\to M$.
Since $\mu$ is also proper the quotient $M = X/\Omega_s G$
is a finite-dimensional compact manifold \cite{AMM}.
The action of $L_sG$ on $X$ induces an action
of $L_s G/\Omega_s G = G$ on $M$.  Furthermore, if we identify the space
$\dlg$ with a space of connections on the trivial principal
$G$-bundle on the circle, the map
${\rm Hol}\colon \dlg \to G$ given by the value at time $1$ of
the holonomy of the connection is a fibration with fibre $\Omega_s G$.
Using this fibration, the $L_sG$--equivariant moment map
$\mu \colon X \to \dlg$ induces a $G$--equivariant map
$\Phi \colon M \to G$, which intertwines the $G$ action on
$M$ with the adjoint $G$ action on $G$.  We thus obtain
a commuting square

\begin{equation}\labell{sqr}
\vcenter{\xymatrix{
X
\ar[r]^{\mu}\ar[d]^{\pi}&
\dlg
\ar[d]^{\rm Hol} &
\\
M
\ar[r]^{\Phi} &
G&
}}
\end{equation}

\noindent showing that the principal $\Omega_s G$-bundle $\pi \colon X\to M$ is
the pullback under $\Phi$ of the contractible principal $\Omega_s G$-bundle
${\rm Hol} \colon \dlg \to G$.

Alexe'ev, Malkin and Meinrenken \cite{AMM} give a set of
conditions a $G$-space $M$ must satisfy in order to arise
from a Hamiltonian $L_sG$ space by this construction.
Given a compact Lie group $G$, choose an invariant
inner product $\langle \cdot , \cdot \rangle$ on $\fg$. 
Let $\theta_L$ and $\theta_R$
denote the left and right invariant Maurer-Cartan
forms,  respectively. and let $\xi ={\frac {1}{12}}\langle\theta_L,[\theta_L,\theta_L]\rangle
={\frac {1}{12}}\langle\theta_R,[\theta_R,\theta_R]\rangle \in \Omega^3(G)$ denote
the bi-invariant $3$--form; if $G$ is simple, $\xi$
generates $H^3(G)$.  We consider $G$ as a $G$-space by equipping
it with the adjoint action.
A {\bf quasi-Hamiltonian $G$-space} (or {\em q-Hamiltonian $G$-space})
$(M,\sigma,\Phi)$
is a compact $G$--manifold $M$, along with a $G$-equivariant
map $\Phi \colon M \to G$, and a two-form $\sigma \in \Omega^2(M)$,
satisfying

\begin{enumerate}
\item $\Phi^*\xi = d \sigma$
\item $i_{\eta^M} \sigma = \langle {\frac 1 2} \Phi^* (\theta_L + \theta_R),\eta\rangle$ for all $\eta \in \fg.$
\item ${\rm ker} \sigma|_p= \{ \eta^M_p \mid \eta \in {\rm ker}(Ad_{\Phi(p)}+ 1)\}$ for all $p  \in M.$
\end{enumerate}

According to \cite{AMM}, Theorem 8.3, there is a one-to-one correspondence
between Hamiltonian $L_sG$-spaces $(X,\omega,\mu)$ with proper
moment map and quasi-Hamiltonian $G$-spaces $(M,\sigma,\Phi)$;
the quasi-Hamiltonian $G$-space $M$ associated to a Hamiltonian
$G$-space $X$ is given by $X/\Omega_sG$, the moment map $\Phi$ is
the map induced on $M$ by $\mu$, and the two-form $\sigma$
is given by formula (31) in \cite{AMM}.
Given a Hamiltonian $L_sG$-space
$X$ with proper moment map, the symplectic quotient
$X \mod L_sG$ is given in
terms of the corresponding quasi-Hamiltonian
$G$-space $M$ by $X \mod L_s G = M \mod G:=\Phi^{-1}(e)/G$.

\mute{
Let $G$ be a compact Lie group.
Choose an invariant
inner product $\langle,\rangle$ on $\fg$.  Let $\theta_L$ and $\theta_R$
denote the left and right invariant Maurer-Cartan
forms,  respectively.
Let $\xi =\langle\theta_L,[\theta_L,\theta_L]\rangle
=\langle\theta_R,[\theta_R,\theta_R]\rangle \in \Omega^3(G)$ denote
the bi-invariant $3$--form; if $G$ is simple, $\xi$
generates $H^3(G)$.  
A {\bf quasi-Hamiltonian G-space} (or {\bf quasi-Hamiltonian $G$-space})
is a triple $(M,\sigma, \Phi)$ where $M$ is a compact $G$--manifold, 
$\sigma \in \Omega^2(M)$ is a two-form on $M$, and $\Phi \colon M \to G$ is
a $G$-equivariant map, which satisfy
\begin{enumerate}
\item $d \sigma = \Phi^*\xi$
\item $i_{\xi^M} \sigma = {\frac 1 2} \Phi^* (\theta_L + \theta_R).$
\item ${\rm ker} \sigma|_p= \{ \eta^M_p \mid \eta \in {\rm ker}(Ad_{\Phi(p)}+ 1)\}
\ \forall p  \in M$
\end{enumerate}

Let $(X,\omega,\mu)$ be a Hamiltonian $L_s G$ space with proper
moment map.
The normal subgroup $\Omega_s G \subset L_sG$ given
by the based maps in $L_s G$
$$\Omega_s G = \{\lambda \in L_s G  \mid \lambda(0) = e\}$$
acts freely on $L_s {\mathfrak g}^*$.
Since $\mu$ is $L_s G$ equivariant,   $\Omega_s G$ acts freely on $X$;
since $\mu$ is also proper the quotient $M = X/\Omega_s G$
is a finite-dimensional compact manifold \cite{AMM}.
The action of $L_sG$ on $X$ induces an action
of $L_s G/\Omega_s G = G$ on $M$.  Furthermore, if we identify the space
$\dlg$ with a space of connections on the trivial principal
$G$-bundle on the circle, the resulting time-one holonomy map
${\rm Hol} \colon \dlg \to G$ is a fibration with fibre $\Omega_s G$.
Using this fibration, the $L_sG$--equivariant moment map
$\mu \colon  X \to \dlg$ induces a  map
$\Phi \colon M \to G$, which  is equivariant with respect to
the adjoint action of $G$ on itself.
We obtain a commuting square
\begin{equation}\labell{sqr1}
\vcenter{\xymatrix{
X
\ar[r]^{\mu}\ar[d]^\pi&
\dlg
\ar[d]^{\rm Hol} &
\\
M
\ar[r]^{\Phi} &
G&
}}
\end{equation}
\noindent showing that the principal $\Omega_s G$-bundle ${\rm Hol} \colon X\to M$ is
the pullback under $\Phi$ of the contractible principal $\Omega_s G$-bundle
${\rm Hol} \colon \dlg \to G$.

In \cite{AMM}, they construct a two form $\sigma \in \Omega^2(M)$
and show that $(M,\sigma,\Phi)$ is a quasi-Hamiltonian $G$-space.
Conversely, they show that every quasi-Hamiltonian $G$-space
arises from  a Hamiltonian $L_sG$ space with a proper moment map
in the fashion described above.

If $e$ is a regular value of $\Phi$, then
the {\bf symplectic quotient} $M \mod G = \Phi\inv(e)/G$
is naturally a symplectic orbifold.  Again,
$G$ acts locally freely on $\Phi\inv(e)$ and  the
equivariant cohomology ring $H^*_G(\Phi\inv(e))$ coincides
with the cohomology ring $H^*(M \mod G)$.
The resulting map $\kappa \colon H_B^*(M) \to H^*(M \mod G)$
is the quasi-Hamiltonian analog of the Kirwan map.
Moreover,  if $(M,\sigma,\Phi)$ is the quasi-Hamiltonian
$G$-space associated to $(X,\omega,\mu)$,
then  $M \mod G$ is
naturally isomorphic to the reduced space $X \mod L_sG$.
}

\begin{Example}
The quasi-Hamiltonian $G$-space  corresponding to $X= LG/G$ in
Example \ref{Ex 1} is a point.
The quasi-Hamiltonian $G$-space corresponding to $X = LG/T$ in
Example \ref{Ex 2} is the
conjugacy class $C_{\exp \xi}$. 
In both cases, the $G$-valued moment map
is given by inclusion.
The quasi-Hamiltonian $G$-space corresponding to Example \ref{Ex 3} is
$M=G^{2h}$, equipped
with the moment map $\Phi \colon M\to G$ given by
$\Phi(a_1,\dots,a_h,b_1,\dots,b_h) = \prod_{i=1}^h [a_i,b_i]$
(see \cite{AMM}).
\end{Example}

We are now ready to state a version of  our main
theorem for quasi-Hamiltonian $G$-spaces.
For any $G$-space $X$, let $C^*_G(X)= C^*(X_G)= C^*(X\times_G EG)$ denote
the singular cochain complex. 
Consider the fibration $p \colon G \times_G EG \to BG$,
and let $j \colon BG \to G \times_G EG$ denote the inclusion
of $BG$ as $\{e\} \times_G EG$.
Because  the action of $G$ on itself is equivariantly formal (see
\cite{gs2}, p. 186), we can find closed co-chains $b_i \in C_G^*(G)$
whose restrictions to the fiber $G$  generate
the cohomology of $G$ as a ring.
By replacing each $b_i$, if necessary, by $b_i - p^*(j^*(b_i))$,
we may assume that $j^*(b_i) = 0$.
If there exist $a_i \in C_G^*(M)$ so that 
$\Phi^*(b_i) = d a_i$
then $a_i|_{\Phi\inv(e)}$ is closed.

\begin{Theorem} [Surjectivity for quasi-Hamiltonian $G$-spaces]
Let $G$ be a compact simply connected Lie group.
Let $b_i \in C^*_G(G)$ satisfy $j^*(b_i) = 0$ and
generate the cohomology of the fiber $G$ of the fibration
$p$ as a ring.

Let $(M,\sigma,\Phi)$ be a quasi-Hamiltonian $G$-space.
Assume that there exist co-chains $a_i \in C^*_G(M)$
such that $d a_i = \Phi^*(b_i)$.

Then $H_G^*(\Phi\inv(e))$ is generated as a ring
by the image of the restriction $H_G^*(M) \to H_G^*(\Phi\inv(e))$
along with the classes $[a_i|_{\Phi\inv(e)}]$.

In particular, if $e$ is a regular value of $\Phi$,
then $H^*(M \mod G)$ is generated as a ring
by the image of the Kirwan map $\kappa$ 
along with the classes $[a_i|_{\Phi\inv(e)}]$.

\end{Theorem}
\begin{Remark} The additional assumption made in our theorem,
that the pull-back $\Phi^*(b_i)$
is  exact for all $i$, is satisfied for
the quasi-Hamiltonian $G$-spaces given in Examples \ref{Ex 1}
,\ref{Ex 2}, and \ref{Ex 3}.
We do not know of
a proof that this condition follows in general
from the definition 
of a quasi-Hamiltonian $G$-space.  However, in the case
where $G=SU(2),$ this is the case, and we have the following
corollary.
\end{Remark}
\begin{corollary}\labell{su2case}
Let $(M,\sigma,\Phi)$ be a quasi-Hamiltonian $SU(2)$-space.

Then $H_{SU(2)}^*(\Phi\inv(e))$ is generated as a ring
by the image of the restriction $H_{SU(2)}^*(M) \to H_{SU(2)}^*(\Phi\inv(e))$
along with the class $[\sigma|_{\Phi\inv(e)}]$.
In particular, if $e$ is a regular value of $\Phi$,
then $H^*(M \mod G)$ is generated as a ring by 
the image of $\kappa$ and the class of the symplectic form.
\end{corollary}

\begin{proof}
In the Cartan model, the equivariant cohomology of a $G$-manifold $M$ is
given by the cohomology of the
complex
$((\Omega^*(M) \otimes S(\fg^*))^G, \delta=d + i_{\eta^M})$,
where $S(\fg^*)$ is the algebra of symmetric polynomials on $\fg$, and
$\eta^M$ denotes the vector field on $M$ corresponding
to $\eta$ for each $\eta \in \fg$; 
see e.g. \cite{gs2}.
If $G$ is a simple group, the generator of $H_G^3(G)$ is
given by the differential form $\xi + {\frac 1 2} \Phi^* (\theta_L +
\theta_R).$
Conditions (1) and (2) in the definition of a quasi-Hamiltonian
$G$-space above state that $\delta\sigma =  
\Phi^*(\xi + {\frac 1 2} \Phi^* (\theta_L + \theta_R)).$
We now apply Theorem 3.

\end{proof}

\begin{proof}[Proof of Theorem 3]
Consider the fibration $q \colon  L_s{\mathfrak g}^* \times_G EG \to BG$,
and let $i \colon BG \to L_s{\mathfrak g}^* \times_G EG$ denote the inclusion
of $BG$ as $\{0\} \times_G EG$.
Since ${\rm Hol} \circ i = j$,
$i^*({\rm Hol}^*(b_i)) = j^*(b_i) = 0 $.
Since $i$ is
a homotopy equivalence,
there exist $c_i \in C_G^*(L_s{\mathfrak g}^*)$ such
that $d c_i = {\rm Hol}^*(b_i)$. 
By replacing each $c_i$, if necessary, by 
$c_i - q^*(i^*(c_i))$, we may assume that $i^*(c_i) = 0$.

Fix a point $x \in BG$, and
let $\gamma_i \in C^*(L_s \fg^*)$ and $\beta_i \in C^*(G)$ 
denote the restrictions of $c_i$ and $b_i$ to a fiber
of $q\inv(x)$ and $p\inv(x)$, respectively.
We now consider the principal $\Omega_s G$ bundle
${\rm Hol} \colon L_s \fg^* \to G$.
Then $d \gamma_i = {\rm Hol}^*(\beta_i)$.
Since the space $L_s{\mathfrak g}^*$ is  contractible
and $H^*(G)$ is generated by the cohomology classes $\beta_i$,
which are classes of odd degree, a spectral sequence
argument shows that 
the restriction of  $\gamma_i$ to each
fiber generates the cohomology $H^*(\Omega_s G)$ as a ring.
Hence, the restriction of $c_i$ to each fiber 
generates the cohomology $H^*(\Omega_s G)$ as a ring.

Now note that $e_i = \mu^*(c_i) -  \pi^*(a_i)$ are closed
cochains on $X \times_G EG$ because $d(e_i) = \mu^*(d c_i) - \pi^*(d a_i)
= \mu^*({\rm Hol}^*(b_i)) - \pi^*(\Phi^*(b_i)) = 0$.
Moreover, the restrictions of the $e_i$ to each fiber of the
fibration $X \times_G EG \to M \times_G EG$
generate the cohomology of the fibre $\Omega_sG$ as a ring.
By the Leray-Hirsch theorem,
the cohomology $H_G^*(X)$ is then generated as a
ring by the $e_i$ and the pull-backs of classes on $H^*_G(M)$
under $\pi$.
Therefore, the image of the restriction map $H_G^*(X) \to H_G^*(\mu\inv(0))$
is also generated as a ring by the images of these classes.
This image is just the image of the restriction
map $H_G^*(M) \to H_G^*(\Phi\inv(e))$,
along with the restrictions $e_i|_{\Phi^{-1}(e)}$ (recall that
$\Phi^{-1}(e) = \mu^{-1}(0)$).
To identify those restrictions, note that
$\pi^*(a_i)|_{\mu^{-1}(0)} = a_i|_{\Phi^{-1}(e)}$.
On the other hand, $\mu^*(c_j)|_{\mu^{-1}(0)}=0$
since the restriction $i^*c_j$ vanishes. Thus
$e_i|_{\mu^{-1}(0)} = a_i|_{\Phi^{-1}(e)}$.

By Theorem 2, the restriction map
$H_G^*(X) \to H_G^*(\mu^{-1}(0))$ is surjective.  It follows that
$H_G^*(\Phi^{-1}(e))$ is generated as a ring by the image
of $H_G^*(M)$ under the inclusion map $k \colon \Phi^{-1}(e) \to M$, along
with the classes $a_i|_{\Phi^{-1}(e)}$.
\end{proof}

\section{Path Spaces}\labell{pathspaces}

In this section we show  that in order to prove
surjectivity for a Hamiltonian $L_sG$-space $X$,
it is enough to prove surjectivity for an associated
infinite-dimensional space $\hX$.
The spaces $X$ and $\hX$ are related in the same way as the
space of paths of Sobolev class  $s >  \frac{1}{2}$ is
related to the space of piecewise smooth paths.
The rest of the paper will be devoted to proving surjectivity for
the space $\hX$.

We first construct this space $\hX$.

Let $\pgeh$
be the space of piecewise smooth based paths on $G$;
that is, piecewise smooth maps $\lambda \colon  [0,1] \to G$ so that
$\lambda(0) = e$,
where  $e \in G$ is the identity.
We define a metric $\delta$ on $\pgeh$, following \cite{Bott}:
given paths $\lambda,\lambda'\in \pgeh$,
$$\delta(\lambda,\lambda') =
{\rm max}_{t \in [0,1]}~d(\lambda(t),\lambda(t'))
+ |J(\lambda)-J(\lambda')|,$$
where $d$ is an invariant metric on $G$ and
$J$ is the length function on $\pgeh$.
The group $G$ acts on  $\pgeh$  by
$(g  \cdot \lambda )(t) = g \inv  \lambda(t) g.$
Let $\hat{\rho} \colon \pgeh \to G$  be the endpoint map,

Let $(X, \omega, \mu)$ be a Hamiltonian $L_s G$-space
with a proper moment map.
Let $(M,\sigma,\Phi)$ be the associated quasi-Hamiltonian $G$--space.
Note that $X$ can be reconstructed from $\Phi \colon M \to G$ since
$$ X = \{(\gamma,m ) \in  L_s\fg^* \times M  \mid \Phi(m) = \rm{Hol}(\gamma)\}.$$
We introduce a new space $\hX$, which we call the
{\bf infinite approximating space} associated to the Hamiltonian-$L_s G$ space
$(X,\omega, \mu)$,
given by
$$\hX :=\{(\lambda,m) \in  \pgeh \times M  \mid  \Phi(m) = \hat{\rho}(\lambda)\}.$$
The space $\hX$  comes equipped with the diagonal  $G$ action. 
Define the  {\bf energy  function} $\hf \colon  \hX \to \R$  by
$$\hf(\lambda,m) = \int_{[0,1]}
\left|\lambda^{-1}{\frac{d\lambda}{dt}}\right|^2 dt.$$

The main result of this section is the following.

\begin{Proposition}\labell{approxi}
Let $G$ be a compact Lie group, and $L_s G$ be the corresponding
loop group, where $s > 1/2$.
Let $(X,\omega,\mu)$ be a Hamiltonian $L_s G$-space with
proper moment map.
Let $\hX$ be the associated infinite approximating
space, and let $\hf \colon \hX \to \R$ be the energy function.

The restriction map
$H^*_G(X) \to H^*_G(f\inv(0))$ is surjective
if and only if the  restriction map
$ H^*_G(\hX) \to H^*_G(\hf\inv(0))$ is surjective.
\end{Proposition}

Let $P_eG^*$ be the space of continuous maps from $[0,1]$ to $G$,
endowed with the uniform topology.
The group $G$ acts on  $P_eG^*$  by
$(g  \cdot \lambda )(t) = g \inv  \lambda(t) g.$
Let ${\rho}^* \colon P_eG^* \to G$  be the endpoint map.
Define the space $X^*$
$$X^* =\{(\gamma,m) \in  P_eG^* \times M  \mid \Phi(m) = \rho^*(\gamma)\}.$$
Define the $G$ action on $X^*$ as in the case of $\hX$. 

Note that the  spaces  $\pgeh$ and $P_eG^*$ are
both  equivariantly contractible.
The space $\hX$ comes equipped with   an equivariant map
$\hmu  \colon \hX\to \pgeh$  defined by
$\hmu(\lambda,m) = \lambda$.
Similarly, $P_eG^*$ comes equipped with  
an equivariant map $\mu^* \colon X^* \to P_eG^*$.
Let  ${\bf e}$ denote the constant path at the
identity in both $\pgeh$ and $P_eG^*$

There is a natural equivariant inclusion map $\hat{i} \colon \hX \to X^*$.
Similarly, where $s > 1/2$,
the map $\hol \colon L_s\fg^*\to P_eG^*$, which sends
every $\fg-$valued one form to its holonomy,
induces a equivariant map $i \colon X \to X^*$.

To prove  Proposition \ref{approxi},  
note that
$f^{-1}(0)=\mu^{-1}(0) = {\mu^*}\inv({\bf e})=\hmu\inv({\bf e})= \hf^{-1}(0)$.
Moreover, the maps $i \colon X \to X^*$ and $\hat{i} \colon \hX \to X^*$ commute
with the restrictions to these subspaces.
Therefore, it is enough to prove the following result.

\begin{Lemma}
The maps $i \colon X \to X^*$ and $\hat{i} \colon \hX \to X^*$ induce isomorphisms in
equivariant cohomology.
\end{Lemma}

Let $\ogh$ be the group of piecewise smooth based loops on $G$, i.e.,
$$\ogh= \{ \lambda \in \pgeh \mid  \lambda(0) = \lambda(1)  = e\}.$$
Let $\ogs$ be  the space of continuous based loops on $G$.
The endpoint maps $\hat{\rho} \colon  \pgeh \to G$ and $\rho^*  \colon \pges \to G$
are (Serre) fibrations, with fibers  $\ogh$ and $\ogs$, respectively.
The natural inclusion map $j \colon \pgh \to P_eG^*$ is a $G$-equivariant
map which preserves these fibrations.

Considering now the exact homotopy sequence of the
fibrations, we get a map of long exact sequences:

\begin{equation}\label{CD}
\vcenter{\xymatrix{
\dots
\ar[r] &
\pi_n(\ogh)
\ar[r] \ar[d] &
\pi_n(\pgeh)
\ar[r] \ar[d]^{j_*} &
\pi_n(G)
\ar[r]  \ar[d]^{=} &
\dots
\\
\dots
\ar[r] &
\pi_n(\Omega G^*)
\ar[r] &
\pi_n(P_eG^*)
\ar[r] &
\pi_n(G)
\ar[r] &
\dots
\\
}}
\end{equation}

\noindent The map on the right is the identity, and the map in the middle
is an isomorphism because both $P_eG^*$ and $\pgh$ are
contractible.  By the five lemma, the inclusion $\ogh\to \Omega G^*$ also
induces an isomorphism in all homotopy groups.

There are
natural fibrations $\hX \to M$ and $X^* \to M$ with fibers $\ogh$ and
$\ogs$, respectively.
The inclusion  $\pgh \to P_eG^*$ is $G$-equivariant, and
preserves these fibrations.
Hence, we  obtain a  map of the homotopy
quotient fibrations from
$\hX \times_G EG \to M \times_G EG$ to $X^* \times_G EG
\to M \times_G EG$.

We  then obtain a map of long exact sequences

\begin{equation}\label{CDE}
\vcenter{\xymatrix{
\dots
\ar[r] &
\pi_n(\ogh)
\ar[r] \ar[d] &
\pi_n(\hX\times_G EG)
\ar[r] \ar[d]&
\pi_n(M\times_G EG)
\ar[r]  \ar[d] &
\dots
\\
\dots
\ar[r] &
\pi_n(\Omega G^*)
\ar[r] &
\pi_n(X^*\times_G EG)
\ar[r] &
\pi_n(M\times_G EG)
\ar[r] &
\dots
\\
}}
\end{equation}

The map on the right side of this sequence is the identity; the map on the
left is an isomorphism for all $n$ by our previous argument.  Thus the map
$i$ induces an isomorphism in all homotopy groups.  The
Whitehead theorem shows that the inclusion $\hX \to X^*$ induces an
isomorphism in equivariant cohomology.

A similar argument shows
that the inclusion $X \to X^*$ also induces an
isomorphism in equivariant cohomology.

\section{Reduction to finite dimensions}

The previous section showed that surjectivity for a Hamiltonian
$L_s G$-space $X$ follows from surjectivity for the associated
infinite approximating space $\hX$.
In this section, we show that the infinite approximating
space can itself be approximated by a sequence of
finite-dimensional manifolds $Y_n$.
We show that surjectivity holds
for $\hX$ if it holds for each $Y_n$.

Let $G$ be a compact Lie group.
Let $(X,\omega,\mu)$ be a Hamiltonian $L_s G$ space with proper moment
map.  Let $(M,\sigma,\Phi)$ be  the associated quasi-Hamiltonian $G$-space.
For each positive integer $n$, let
$$X_n  \colon= \left\{ (g_1,\ldots,g_n, m) \in G^{n} \times M \mid g_n =
\Phi(m)\right\}.$$
Note that $X_n$ is diffeomorphic to $G^{n-1} \times M$.
The space $X_n$ comes equipped with the diagonal $G$ action,
where $G$ acts on each copy of $G$ by the adjoint action.
The {\bf  energy function}
$f_n \colon  X_n \to \R$ is  given by
$$f_n(g_1,\ldots,g_{n},m)  =
n  \rho(e,g_1)^2 + n \rho(g_1,g_2)^2 + \cdots
+ n \rho(g_{n-1},g_{n})^2.  $$
Note that the energy function $f_n$ is $G$-invariant.
There exists a positive number
$\orho$ such that any two points $p,q \in G$ of distance
$d(p,q) < \orho$ may be connected by a unique shortest
geodesic.
We define the {\bf finite approximating manifold} $Y_n$
associated to $(X,\omega,\mu)$ by
$$Y_n  :=  f_n\inv(-\infty,  \frac{1}{2} n \orho^2).$$

We can now state the main result of this section:

\begin{Proposition}\labell{finite}
Let $(X,\omega,\mu)$ be a Hamiltonian $L_s G$ space with
proper moment map.
Let $\hX$ be the associated infinite approximating
space, and let $\hf \colon \hX \to \R$ be the energy function.
Let $Y_n$ be the finite approximating manifolds,
and let $f_n \colon Y_n \to \R$ be the energy functions.

If the restriction maps $H_G^*(Y_n)\to H_G^*(f_n \inv(0))$ are all
surjections, then the restriction map  $H_G^*(\hX) \to H_G^*(\hf
\inv(0))$
is a surjection.
\end{Proposition}

Given $p=(g_1,\ldots,g_n,m) \in Y_n$, let $g_0 = e$.
Then $$f_n(g_1,\ldots,g_n,m) =
n \left( \sum_{i=0}^{n-1} \rho(g_{i},g_{i+1})^2 \right) < \frac{1}{2} n
\orho^2.$$
Therefore,  $\rho(g_i,g_{i+1}) \leq \sqrt
{\frac
{\frac{1}{2} n \orho^2}
{n}} <
\orho$ for all $0 \leq i < n$.
Hence, there is a unique  minimal
geodesic  between $g_i$ and $g_{i+1}$.
These geodesics can be joined  to form a path
$\lambda_p$ such that $\lambda_p(\frac{i}{n}) = g_n$.
Define $\beta\colon Y_n \to  \hX$
by
$$\beta(p=(g_1,\ldots,g_n,m)) =  (\lambda_p,m).$$
It is clear from the definition that
$$\hf  \circ \beta = f_n.$$
Hence, $\beta$ is a map from $Y_n$ to
$\hf\inv(-\infty, \frac{1}{2} n \orho^2) \subset \hX$,  
where
$$\hf\inv(-\infty, \frac{1}{2} n \orho^2) =  
\{ x \in \hX \mid \hf(x) < \frac{1}{2}n \orho^2 \}. $$
Propostion \ref{finite} now follows immediately
from the lemma below, which
is adapted from \cite{Bott} and \cite{Milnor}.

\begin{Lemma}\labell{approx}
For any natural number  $n$,
$\beta \colon Y_n \to \hf^{-1}(-\infty, \frac{1}{2}n \orho^2)$
is a  $G$-equivariant homotopy equivalence.
\end{Lemma}

\begin{proof}
First, we define $\alpha \colon
\hf\inv(-\infty, \frac{1}{2} n \orho^2) \to Y_n.$
Given $(\lambda, m) \in \hX$, define  $\alpha(\lambda,m) \in X_n$ by
$$\alpha(\lambda,m) =
\left( \lambda(1/n),\lambda(2/n), \ldots,
\lambda(1), m\right).$$
This function is well-defined,
since, by assumption $\lambda(1) = \Phi(m)$.
It is also  continuous.
Moreover,
this map   takes 
$\hf\inv(-\infty, \frac{1}{2} n \orho^2)$
to $Y_n$
since
$$\int_{\left[\frac{i}{n},\frac{i+1}{n}\right]}
\left|{\frac{d\lambda}{dt}}\right|^2 dt
\geq n \rho
\left(\lambda\left(\frac{i}{n}\right),\lambda \left(\frac{i+1}{n}
\right) \right)^2$$
for all $0 \leq i < n$.

As in \cite{Bott}, one may
construct a homotopy $D_t$ from 
$\hf\inv(-\infty, \frac{1}{2} n \orho^2) \subset \hX$
to itself so
that $D_0$ is the identity, and $D_1 = \beta \circ \alpha$.
To construct this, one simply deforms the segment of
$\lambda$ between $\lambda(\frac{i}{n})$ and $\lambda(\frac{i+1}{n})$
into the shortest geodesic joining
$\lambda(\frac{i}{n})$ and $\lambda(\frac{i+1}{n}).$
The intermediate paths will be  the shortest geodesic joining
$\lambda(\frac{i}{n})$ and $\lambda(\frac{i}{n} + \epsilon),$
followed by the original curve between
$\lambda(\frac{i}{n} + \epsilon)$ and $\lambda(\frac{i+1}{n}).$

Finally, since  $\alpha \circ \beta$ is already the identity,
we are done.
\end{proof}

\section{Properties of the finite-dimensional approximation}

To prove that the energy function $f_n \colon Y_n \to \R$ on the
finite approximating manifold is Morse in the sense of Kirwan,
one needs 
to find the critical sets of $f_n$ and to
compute the Hessian on each critical set.
In this section, we begin these computations.

The results in this section do not depend on the details
of the properties
of quasi-Hamiltonian $G$-spaces: Given any compact, connected
$G$-manifold $M$, along with an equivariant
map $\Phi \colon M \to G$,
define $X_n$, $f_n$, and $Y_n$ as in Section 4 above.  The results
in this section will then hold.

\begin{Lemma}\labell{critset}
Consider $y = (g_1,\dots,g_{n}, m) \in Y_n$.
Let $\beta(y) = (\lambda,m) \in \hX.$
The function $f_n$ is critical at $y$ if and only if
following conditions hold:
\begin{enumerate}
\item The path $\lambda$ is a geodesic. 
\item The velocity of $\lambda$ at time $1$
is perpendicular to the image of $\Phi_* \colon T_m M \to T_{\Phi(m)} G$.
\end{enumerate}
\end{Lemma}

\begin{proof}
(Compare with \cite{Bott}, p.\ 319.)
Let $X = (X_1,\ldots, X_{n-1}, X_n, Y) \in T_y Y_n$ be a tangent vector.
Note that  $X_n = \Phi_*(Y)$, but otherwise
the $X_i$  are independent.

Let $g_0 = e$ and let $X_0  \in T_e(G)$ be the zero vector.
Let $s_i$ denote the unique shortest geodesic between $g_i$ and
$g_{i+1}$ for all $0 \leq i < n$, parametrized
so that $s_i(0) = g_i$ and $s_i(1) = g_{i+1}$.
Let $\dot{s}_i(t)$ denote the unit tangent vector of $s_i$ at
time $t$.  Consider a one-parameter family of paths
$\theta(s,t)$ so that $\theta(0,t) = s(t)$ and
$\frac{\partial \theta }{\partial s}(0,t) = W(t),$ where $W(t)$
is the unique Jacobi field which is continuous
on $s(t)$ and smooth except possibly at $\frac{i}{n}$,
such that $W(i/n) = X_i$.  

By the first variation formula
$$X ( \rho^2(g_i,g_{i+1})) =
2  \left< \dot{s}_i(1),  X_{i+1} \right>
- 2  \left< \dot{s}_i(0),  X_i \right>.$$
Summing the terms, we get
$$\frac{1}{2n} X ( f_n) =   \sum_{i = 0}^{n-2}
\left< \dot{s}_i(1) -  \dot{s}_{i+1}(0) , X_{i+1} \right>
+  \left<  \dot{s}_{n-1}(1) , \Phi_*(Y) \right>.$$

This is zero for all $X \in T_yY$ exactly if
$\dot{s}_i(1) =  \dot{s}_{i+1}(0)$ for all $i$
and  $\Phi_*(T_y Y_n)$ is perpendicular to
$ \dot{s}_{n-1}(1)$
\end{proof}

\begin{Definition}\labell{hessian}
Consider any $m \in M$ and $\xi \in T_{\Phi(m)} G$ such that $\xi$ is
perpendicular to the image  $\Phi_*(T_m M)$.
There exists a symmetric bilinear form
$$H ^\xi \Phi \colon T_m M \times T_m M \to \R$$
defined as follows:
Given $Y$ and $Y'$ in $T_m M$, choose a smooth map $\alpha$ from a
neighborhood of $(0,0)$ in $\R^2$ to $M$ so that
$$ \alpha(0) = m,  \ \   \frac{\partial \alpha}{\partial s}(0,0) = Y,
 \ \  \mbox{and} \ \
 \frac{\partial \alpha}{\partial t}(0,0) = Y',$$
where $s$ and $t$ are the coordinates on $\R^2$.
Then $$H^\xi \Phi(Y,Y') = \
\left< \frac{D}{D s} \frac{\partial}{\partial t} ( \Phi \circ
\alpha) (0,0), \xi \right>.$$

Here we have denoted by $D/Ds$ the covariant derivative
given by the Levi-Civita connection associated
to the bi-invariant metric on $G$.
\end{Definition}

The proof that this is well-defined is analogous to
the proof that the Hessian of a function is well-defined at a critical
point.

\begin{Lemma}
Let  $y = (g_1,\ldots,g_n,m) \in Y_n$ be a critical point of $f_n$
and let $\beta(y) = (\lambda,m)$.
Consider $\eta = (X_1,\ldots,X_n,Y)$ and
$\eta' = (X'_1,\ldots,X'_n,Y') \in  T_y Y_n$. Let  $W$ be the unique
Jacobi field on  $\lambda$ which is smooth except possibly at
$g_i$, and such that $W(0) =  0$ and $W(i/n) = X_i$ for all $1 < i \leq n$.
Let $\Delta_i \frac{D W}{D t}$ denote
the difference between the left and right hand limit of
$\frac{D W}{D t}$ at $g_i$.
The Hessian $H f_n$ of $f_n$ at $y$ is given by the following formula:
\begin{equation}\labell{Hessianfn}
H f_n (\eta, \eta') =
\sum_{i < n} \left<  \Delta_i \frac{DW}{Dt}, X'_i \right>
+ \left<  \frac{D W}{D t}(1), X'_n \right>
+ H^\xi \Phi(Y,Y'). \end{equation}
\end{Lemma}

\begin{proof}
The second variation formula for fixed endpoints is
easily extended to the case where
variations in the path are allowed at $t=1$.  Given
a geodesic $\lambda(t)=\exp(\xi t)$
 on the group $G$, and let
$\theta(s,t)$ be a variation of $\lambda$; that is
$\theta(0,t)=\lambda(t),$ and $\theta(s,0)=e$ for all
$s$.  Suppose $\theta$ is continuous everywhere and
smooth except possibly where $t=i/n, i=1,\cdots,n-1,$
Write 
$V=\theta_*(\partial/\partial t),$
$Y=\theta_*(\partial/\partial s),$ and $W=Y|_{s=0}$.  The
vector field $W(t)$
is a broken Jacobi field along $\lambda;$ we write
$X_i = W(i/n).$  Let
$E(s) :=1/2 \int_0^1 |\partial \theta/\partial t|^2 dt.$
The familiar
second variation formula 
(e.g. \cite{Milnor}, Theorem 13.1)
becomes
$d^2E/ds^2|_{s=0} = 
\sum_{i < n} \left<  \Delta_i \frac{DW}{Dt}, X_i \right>
+ \langle D_Y Y, V\rangle|_{s=0,t=1}
+\langle W(1), \frac{DW}{Dt}(1)\rangle|_{s=0,t=1}.$  The result follows.

\end{proof}

Since $f_n$ is a $G$-invariant function,  equation (\ref{Hessianfn})
gives the following result.

\begin{Lemma}\labell{invar}
Let  $y = (g_1,\ldots,g_n,m) \in Y_n$ be a critical point of $f_n$
and let $\beta(y) = (\lambda,m)$.
Given $\gamma\in \fg$,
let $\hat{\gamma} \in T_m M$ be the value at $m$
of the vector field $\gamma^M$ on $M$ induced by $\gamma$.
Let $W$ be the Jacobi field
on the path $\lambda$ induced by conjugation by $exp(s \gamma)$.
Then
$$ \left<DW/Dt(1), \Phi_*(Y') \right> + H^\xi \Phi(\gamma_M,Y')  = 0$$
for all $Y' \in T_m(M)$.
\end{Lemma}

\section{Local normal forms for quasi-Hamiltonian $G$-spaces}

The main goal of this section is to prove 
a local normal form theorem for quasi-Hamiltonian
$G$-spaces (Proposition \ref{qnormalform}).
The proof of this theorem involves combining 
the local normal form theorem for Hamiltonian $G$-spaces \cite{gs}
with results from \cite{AMM} which show that a quasi-Hamiltonian
$G$-space is locally modelled on a Hamiltonian $G$-space.
The normal form theorem of Proposition \ref{qnormalform} 
will allow us to perform explicit calculations of the
Hessians of  
the energy functions on the finite approximating manifolds, which
we will use to show that these functions are 
Morse in the sense of Kirwan.

Let $G$ be a compact Lie group, and let $\langle \cdot, \cdot \rangle$
be an invariant metric on $\fg$. 
Given a subpace $\fh \subset \fg$, let $\fh^\perp \subset \fg$
denote its  metric orthogonal complement.
We will also use the metric to
identify $\fh^*$ with $\fh$.

Let $(M,\sigma,\Phi)$ be a quasi-Hamiltonian $G$-space.
Given a subspace $W \subset T_x M$,  let $W^\sigma \subset T_x M$
denote the subpace of $\sigma$ orthogonal vectors.
The {\bf symplectic slice} at $p \in M$ is the vector space
$$ V = (T_p \cO)^\sigma / (T_p \cO \cap (T_p \cO)^\sigma),$$
where $\cO = G \cdot p$ is the $G$ orbit of $p$.
Since by the axioms for quasi-Hamiltonian $G$-spaces,
the kernel of $\sigma_p$ is contained entirely in $T_p \cO$,
$V$ is a symplectic vector space.
The isotropy representation of $\Stab(p)$ on $T_p M$ induces
a representation on  the symplectic slice, called the {\bf slice
representation.} Our main proposition is the following.

\begin{Proposition}\labell{qnormalform}
Let $(M,\sigma,\Phi)$ be a quasi-Hamiltonian $G$-space.
For any $p \in M$,
let $H= \Stab(p)$,  $K= \Stab(\Phi(p))$, and $V$ be the symplectic
slice at $p$.

There exists a neighbourhood of the orbit $G \cdot p$ which is equivariantly
diffeomorphic to a neighborhood of the orbit $G\cdot [e,0,0]$  in
$$ Y := G \times_H ((\fh^\perp \cap \fk) \times V).$$
In terms of
this diffeomorphism, the $G$-valued moment map
$\Phi \colon M \to G$ may be written as
$$\Phi([g,\gamma,v]) =  \Ad_g (\Phi(p) \exp(\gamma + \phi(v)) ),$$
\noindent where $\phi \colon V \to \fh^*\simeq \fh$ is the moment map for the
slice representation.
\end{Proposition}

To prove Proposition \ref{qnormalform},
we will need the local normal form theorem for Hamiltonian
$G$-spaces
given, for example, in Guillemin-Sternberg \cite{gs},
Section 41.
Given a subalgebra $\fh \subset \fg$, let $\fh^0 \subset \fg^*$
denote its annihilator.

\begin{Proposition}\labell{normalform}
Let $(M,\omega,\mu)$ be a Hamiltonian $G$-space.
For any $p \in M$, let $H= \Stab(p)$, let $K= \Stab(\mu(p))$,
and let $V$ be the symplectic slice at $p$.
There exists a neighbourhood of the orbit
$G \cdot p$ which is equivariantly diffeomorphic to
a neighborhood of
the orbit $G\cdot [e,0,0]$  in
$$ Y := G \times_H ((\fh^0 \cap \fk^*) \times V).$$
In terms of this diffeomorphism, the moment map
$\mu \colon M \to \fg^*$ may be written as
$$\mu([g,\gamma,v]) = \Ad^*_g (\mu(p) + \gamma + \phi(v)),$$
\noindent where $\phi \colon V \to \fh^*$ is the moment map for the
slice representation.
\end{Proposition}

Here, the {\bf symplectic slice} at $p \in M$ is the vector space
$ V = (T_p \cO)^\omega / (T_p \cO \cap (T_p \cO)^\omega),$
where $\cO = G \cdot p$ is the $G$ orbit of $p$.

In order to derive the normal form theorem for quasi-Hamiltonian $G$-spaces
from the normal form theorem for Hamiltonian $G$-spaces, we need the
following 
two theorems about quasi-Hamiltonian $G$-spaces,
both taken from \cite{AMM}, which show that locally, a quasi-Hamiltonian
$G$-space may be modelled on a Hamiltonian $G$-space.

\begin{Proposition}[\cite{AMM}, Remark 3.3]\labell{exp}
Let $(M,\sigma,\Phi)$ be a quasi-Hamiltonian $G$-space.
Let $U \subset {\mathfrak g}$ be a connected neigborhood of $0$
so that the exponential map is a diffeomorphism on $U$,
and let $V= \exp U$.  Then there exists
a Hamiltonian $G$-space $(N,\omega,\nu)$
and an equivariant
diffeomorphism $\psi \colon N\to \Phi^{-1}(V)$,
so that the following diagram commutes

\begin{equation}
\vcenter{\xymatrix{
N
\ar[r]^{\nu}\ar[d]^{\psi}&
{\mathfrak g}^*\simeq{\mathfrak g}
\ar[d]^{\rm exp}&
\\
\Phi^{-1}(V)
\ar[r]^{\Phi|_{{\Phi}^{-1}(V)}}&
G&
}}
\end{equation}

\end{Proposition}

Consider the adjoint action of the Lie group $G$ on itself.
Recall that for any $\zeta \in G$, $U_\zeta \subset K$ is a {\bf
slice for the action of $G$ at $\zeta$} if
$U_\zeta$ is preserved by the action of $K=Z(\zeta)$,
the centralizer of $\zeta$,
and if the natural map 
$G \times_K U_\zeta \to G$ is an equivariant diffeomorphism onto its image.
Note that a slice exists for every $\zeta \in G$.

\begin{Proposition}[\cite{AMM}, Proposition 7.1]\labell{crosssection}
Let $(M,\sigma,\Phi)$ be a quasi-Hamiltonian $G$-space.
Given $\zeta \in G$, let $U_\zeta$ be a slice for the action of $G$
on itself at $\zeta$, and let $Y_\zeta :=\Phi\inv(U_\zeta)$.
Let $K=Z(\zeta)$ be the centralizer of $\zeta$.
The {\bf quasi-Hamiltonian cross-section}
$(Y_\zeta, \sigma|_{Y_\zeta}, \Phi|_{Y_\zeta})$
is a  quasi-Hamiltonian $K$-space.
\end{Proposition}

We can now begin our proof.

\begin{proof}[Proof of Proposition \ref{qnormalform}]
Let $p \in M$ and let $\zeta := \Phi(p)$.
Let $U_\zeta$ be a slice for the action of $G$ on itself at $\zeta$.
Let $Y_\zeta = \Phi\inv(U_\zeta).$
By Proposition \ref{crosssection}, 
$(Y_\zeta, \sigma|_{Y_\zeta}, \Phi|_{Y_\zeta})$
is a  quasi-Hamiltonian $K$-space.
Define
$\Psi \colon Y_\zeta \to K$ by $\Psi(m) = \zeta\inv \Phi(m)$.
Since $\zeta$ is in the center of $K$, the triple
$(Y_\zeta,\sigma|_{Y_\zeta}, \Psi|_{Y_\zeta})$ is
also a quasi-Hamiltonian $K$-space.
Moreover $\Psi(p) = e$.
Therefore, by Proposition \ref{exp},
there exists a
Hamlitonian $K$-space $(N,\omega,\nu)$ 
that is $K$-equivariantly diffeomorphic to a neighborhood of $p$ in $Y_\zeta$,
with the diffeomorphism carrying $\exp(\nu)$ to $\Psi$.
A calculation using the explicit formula given in
\cite{AMM} (equation 31) for the relation between the form $\sigma|_{Y_\zeta}$ and
the form $\omega$ shows that the symplectic slice $V$ at the point corresponding to $p$ 
in $N$ 
is identical with the symplectic slice at $p$ in $Y_\zeta$.
Finally, we apply the local normal form theorem 
for Hamiltonian
$K$-spaces (Proposition \ref{normalform}).
This shows that $N$ is locally equivariantly diffeomorphic
to $K \times_H  ((\fh^0 \cap \fk^*) \times V)$, with the diffeomorphism
carrying the moment map $\nu$ to the map
$[k,\alpha,v] \to k \cdot( \alpha + \phi(v))$, where
$\phi$ is the moment map for the slice representation.
\end{proof}

\section{Properties of quasi-Hamiltonian $G$-spaces}

In this section, we first prove the collection of results immediately
below, which  describe
the first and second derivatives of a $G$-valued moment map.
We conclude by proving a result which will be useful in characterizing
the critical set of the functions $f_n$.  All
these results are proved using the local normal form theorem from
the previous section.

We will use the following notation.  Let $G$ be a compact Lie group,
and let $M$ be a $G$-space.
Given $m\in M$, let $\stab(m):= {\rm Lie}(\Stab(m))$.  Also,
for any $\xi \in \fg$, let $M^\xi := \{m\in M \mid e^{t\xi}\cdot m = m 
\ \forall\ t \in \R\}$.

\begin{Lemma}[\cite{AMM}, Proposition 4.1, part (3)]
\labell{firstorder}
Let $(M,\sigma,\Phi)$ be a quasi-Hamiltonian $G$-space.
Fix  $m \in M $; let $\fh = \stab(m)$.
Then $$ \Phi(m)\inv \Phi_*(T_mM) = \fh^\perp.$$
\end{Lemma}

\begin{Lemma}\labell{hamimage}
Let $(M,\sigma,\Phi)$ be  a quasi-Hamiltonian $G$-space.
Fix $\xi \in  \fg$  and $m \in M^\xi$ such that $\Phi(m) = \exp(\xi)$;
let $\fh = \stab(m)$.
Then $$\Phi(m)\inv \Phi_*(T_m M^\xi) = \fh^\perp \cap \fg^\xi.$$
\end{Lemma}

\begin{Proposition}\labell{commuting}
Let $(M,\sigma,\Phi)$ be  a quasi-Hamiltonian $G$-space.
Fix $\xi \in  \fg$  and $m \in M^\xi$ such that $\Phi(m) = \exp(\xi)$.
If $X\in T_m(M^\xi)$,  then $X$ is in the null-space of $H^\xi\Phi_m$.
\end{Proposition}

\begin{Proposition}\labell{nondeg}
Let $(M,\sigma,\Phi)$ be  a quasi-Hamiltonian $G$-space.
Fix $\xi \in  \fg$  and $m \in M^\xi$ such that $\Phi(m) = \exp(\xi)$.
Given $X \in T_m(M)$  such that $\Phi(m)\inv \Phi_*(X) = 0$,
if $X$ is in  the null-space of the Hessian $H^\xi \Phi_m$, then
$X \in T_m(M^\xi).$
\end{Proposition}

Let $\zeta = \Phi(m)$, $H = \Stab (m)$,  and $K = \Stab (\zeta)$.
Let $V$ be
the symplectic slice at $m$, and let $\phi \colon V \to \fh^*\simeq \fh$
be the moment map
for the slice representation.
By Proposition \ref{qnormalform},  there exists a neighborhood
of $G \cdot m$ equivariantly diffeomorphic to a neighborhood
of $G \cdot [e,0,0]$ in $Y =  G \times_H (\fk \cap \fh^\perp \times V)$.
Under this identification, $\Phi$ is given by
$\Phi([g,\gamma,v])= Ad_g( \zeta \exp (\gamma + \phi(v))).$

A direct computation shows that for
$p = [e,0,0] \in  G \times_H (\fk \cap \fh^\perp \times V)$
and for
any  $(a,x,v) \in \fh^\perp \times \fk \cap \fh^\perp \times V = 
T_p (G \times_H (\fk \cap \fh^\perp \times V))$,
\begin{equation}\labell{derivatives}
\Phi(p)\inv \Phi_* (a,x,v) = (1-Ad_{\zeta})(a) + x.
\end{equation}

Consider the map from $h \colon \fg\to\fg$ given by $h(a) = (1 - \Ad_\zeta)(a)$.
The kernel of $h$ is given by ${\rm ker}(h) = \fk$.
On the other hand, if $b \in \fk$, then $\left< (1 - \Ad_\zeta)(a),b \right>
= \left< a, (1 - \Ad_{\zeta\inv})(b) \right> = 0$ for all $a \in \fg$,
so $h(\fg) \subset \fk^\perp$.
By a dimension count, this implies that
$h(\fg)=\fk^\perp$.
Lemma \ref{firstorder} now follows immediately from equation
(\ref{derivatives}).

We now prove Lemma \ref{hamimage}.
Note  that
$\Phi(p)\inv  \Phi_* (TY^\xi)$ is contained in $\fh^\perp \cap
\fg^\xi$.  On the other hand, since  $\exp(\xi) = \zeta$, 
$\fg^\xi \subset \fk$.  Thus, given $x \in \fh^\perp \cap \fg^\xi$,
we have $(0,x,0) \in T_p Y^\xi$ and $\Phi(p) \inv \Phi_*(0,x,0) = x$.

We now consider the proofs of Proposition \ref{commuting} and 
Proposition \ref{nondeg}.  As these involve the Hessian
$H^\xi\Phi$, we begin with a direct computation of this
Hessian in local coordinates.  To do so, choose $[a,x,v]\in T_{[e,0,0]}Y$
and paths $g(s)$ in $G$, $x(s)$ in $\fh^\perp \cap \fk$ and 
$v(s)$ in $V$ so that $g(0)=e,g'(0) = a, x(0)=0, x'(0) = x, v(0)=0,$ and
$v'(0)=v.$  We write 
$\Phi(s) = {\rm Ad}_{g(s)} \zeta \exp(x(s) + \phi(v(s))$
and compute the quantities 
$H^\xi_L \Phi:=d/ds \langle \Phi(s)^{-1} d/ds (\Phi(s)),\xi\rangle|_{s=0}$
and
$H^\xi_R\Phi :=d/ds \langle (d/ds (\Phi(s)))\Phi(s)^{-1},\xi\rangle|_{s=0}.$  These
quantities correspond to a variant of the definition of the $H^\xi\Phi$
given in Definition \ref{hessian} where the covariant derivative
with respect to the Levi-Civita connection is
replaced by its analog for the connection obtained from the parallelism
arising from left-translation to the origin
and right-invariant translation to the origin, respectively.
It is easy to see that
that $H^\xi_L\Phi = H^\xi_R\Phi$, and since the Levi-Civita
connection is the average of the left- and right-invariant connections,
$H^\xi\Phi = H^\xi_L\Phi = H^\xi_R\Phi$.  
A direct calculation of $H^\xi_L\Phi$ shows that
given  $(a,x,v)$ and $(a',x',v')$ in $T_p Y$, the Hessian of
$\Phi$ at $p$ is given by
\begin{equation} \labell{Hessian}
2  H^\xi \Phi_p( (a,x,v),(a',x',v')) = \left< \phi(v,v'), \xi \right> + \end{equation}
$$ \left<[a, \Ad_\zeta a'] + [a', \Ad_\zeta a] + [x',(1+\Ad_\zeta)(a)]
+ [x, (1 + \Ad_\zeta)(a')] , \xi \right>.$$
Here, $\phi(v,v')$ denotes the $\fh$-valued function
on $V\times V$ which
gives rise to the quadratic moment map $\phi \colon V \to \fh^*\simeq \fh$.

We  now prove Proposition \ref{commuting}.
Given any
$(a,x,v) \in T_p Y^\xi= (\fh^\perp)^\xi \times (\fh^\perp)^\xi \times V^\xi$,
and any $(a',x',v') \in T_p Y$, we may apply
equation (\ref{Hessian}) to obtain
$$ 2 H^\xi \Phi_p ((a,x,v),(a',x',v')) = $$
$$\left<[a, a'] + [a', a] + 2 [x',a] + 2 [x, a'] + \phi(v,v'), \xi \right>= 0,$$
as needed.

We  now  consider Proposition \ref{nondeg}.
By equation (\ref{derivatives}), 
the kernel of $\Phi(p)\inv  \Phi_*$ is
$\fk \cap \fh^\perp \times 0 \times V$.
Furthermore, if $(a,0,v) \in \fk\cap\fh^\perp \times 0 \times V$ is
in the nullspace of $H^\xi\Phi_p$, then
$$2H^\xi\Phi_p ((a,0,v),(0,0,v'))= \langle\phi(v,v'),\xi\rangle=0$$
\noindent for all $v' \in V$.  Thus $v\in V^\xi$ and $(a,0,0)$ is
also in the nullspace of $H^\xi\Phi_p$.  Note that since
$a\in \fk \cap \fh^\perp$ and $\xi \in \fk$,
$[a,\xi] \in \fk \cap \fh^\perp$.
By equation (\ref{Hessian}), 

$$2 H^\xi\Phi_p((a,0,0) , (0,[a,\xi],0)) = 2\langle [a,\xi],[a,\xi]\rangle = 0.$$
\noindent It follows that $[a,\xi]=0$, so that $(a,0,v) \in T_pM^\xi$.

Finally we have the following result.

\begin{Proposition}\labell{goodcriticalsets}
Let $(M,\sigma,\Phi)$ be a quasi-Hamiltonian $G$-space.
Define $\Delta \subset \fg \times M$ by 
$$\Delta:=  \{(\xi,m) \in \fg \times M \mid \xi \in {\rm stab}(m) {\rm and~}
\exp(\xi) = \Phi(m)\}.$$
Fix $r> 0.$
The set $\Delta \cap ((G\cdot \xi) \times M)$ is nonempty for
a finite number of $G$-orbits $G\cdot \xi$ with $|\xi| < r.$
Furthermore, for each $\xi \in \fg$, the set
$\Delta \cap ((G\cdot\xi) \times M)$ has
a finite number of connected components.
\end{Proposition}

\begin{proof}
Since $\Delta$ is closed and $M$ is compact, it will suffice to show that 
for all $(\eta,m) \in \Delta$,
there exists a $G$-invariant neighborhood $U$ of the orbit $G \cdot m$ so that
$\Delta \cap ((G\cdot \xi) \times U)$ is nonempty for only finitely
many orbits $G\cdot \xi$ with $|\xi| < r,$ and 
so that $\Delta \cap ((G \cdot \eta) \times U)$
has finitely many connected components.
So fix $(\eta,m) \in \Delta$.

Let $\zeta := \Phi(m)$, $H:= {\rm Stab}(m)$, and
$K:= Z(\zeta)$.  Since $(\eta,m) \in \Delta$, $\eta \in \fh$.
We can choose a maximal torus $T$ for $G$ so that
$\eta \in \ft$ and $T \cap H$ is a maximal torus for $H$.

Let $V$ be the symplectic slice, and let
$\phi:V \to \fh^* \simeq \fh$ be the moment map for the slice
representation.
By Proposition \ref{qnormalform}, there exists a
$G$-invariant neighborhood $U$ of
$G\cdot m$, an $H$-invariant neighborhood $A$
of $ 0 \in \fh^\perp  \cap \fk$, and
an $H$-invariant
neighborhood $B$ of $0 \in V$ so that $U \simeq G\times_H (A\times B).$
Under this identification, the moment map $\Phi|_U$ is given by
$\Phi([g,a,b]) = {\rm Ad}_g \zeta \exp(a+\phi(b)).$  There exists a $K$
invariant neighborhood $Z$ of $\zeta$ in $K$ so that the natural map
$$ G \times_K Z \to G$$ is an equivariant diffeomorphism onto
its image.
By shrinking $U$
further, if necessary, we may assume that the set $A \times \phi(B)$
is contained in the injectivity radius of the exponential map
$\exp : \fk \to K$ and also that $\zeta \exp(A + \phi(B)) \subset Z$.

We first show that 
$\Delta \cap ((G\cdot \xi) \times U)$ is nonempty for only finitely
many orbits $G\cdot \xi$ with $|\xi| < r.$  
Given a conjugacy class $\Lambda$ of subgroups of $G$, we write
$U_{\Lambda} := \{u \in U\mid {\rm Stab}(u) \in \Lambda\}$.  Because
$U_\Lambda \neq \emptyset$ for only a finite number of conjugacy
classes, it will suffice
to show that for every conjugacy class $\Lambda$,
$\Delta \cap ((G\cdot \xi) \times U_\Lambda)$ is
nonempty for only finitely many orbits $G \cdot \xi$.

Fix $\Lambda$ so that $U_\Lambda \neq \emptyset$.
Choose a representative
$\tilde{H}$ of $\Lambda$ 
so that $T \cap \tilde{H}$ is
a maximal torus of $\tilde{H}$.
Let $B_r(0)$ denote the ball of radius $r$ in $\fg,$ and
define, for $\xi \in \tilde{\fh} \cap \ft \cap B_r(0),$
$$V_\xi:=\{(a,b) \in A \times B \mid  
\Stab(a,b) = \tilde{H} \ \mbox{and} \exp(\xi) = \zeta
\exp(a + \phi(b)) \}.$$ 
The set $\Delta$ is $G$-invariant.  Moreover
every $G$-orbit in $\Delta \cap (\fg \times U_\Lambda)$
contains a point $(\xi,u)$ so
that $u$ is of the form $[e,a,b],$ $\Stab(u) = \tilde{H}$, and
$\xi \in \ft \cap \tilde{\fh}$.
Hence, we will be done if we prove that
$V_\xi$
is nonempty for only finitely many
$\xi \in \tilde{\fh} \cap \ft \cap B_r(0).$ 
let
$\ell \subset \ft$ denote the integral lattice in the torus
$\ft$.
Then, for $\xi \in \tilde{\fh} \cap \ft \cap B_r(0),$
$$V_\xi=\{(a,b) \in \times  A
\times B \mid  \Stab(a,b) = \tilde{H} \mbox{~and}\ 
\eta + a + \phi(b) - \xi \in \ell  \}.$$

Now fix $\lambda \in \ell$.
Since $A \times \phi(B)$ is a bounded subset
of $\tilde{\fh}^\perp$,
the set
$\tilde{\fh} \cap \ft \cap B_r(0) \cap (A \times \phi(B)+ \lambda+\eta)$
contains at most one point, and is empty for sufficiently
large $\lambda$.
This proves that $V_\xi$ is non-empty for only finitely many such $\xi$.

We will now show that $\Delta \cap ((G \cdot \eta \times U)$
has finitely many connected components.
Every orbit in $\Delta \cap ((G \cdot \eta \times U)$ contains
a point $(\xi,u)$ so that $u$ is of the form $[e,a,b]$ and $\xi \in \ft \cap
\fh$.  Since $G \cdot \eta \cap \ft$ is a finite set, it is
enough to check that for all $\xi \in \ft \cap \fh \cap G \cdot \eta$,
the set
$\{ (a,b)  \in A \times B \mid \exp(\xi) = \zeta \exp(a + \phi(b)) \}$
is connected.
Now, $\exp(\xi) = g\inv \zeta g$.
So $\Ad_g(\zeta) = \zeta \exp(a + \phi(b))$.
Hence $\exp(a + \phi(b)) = 1$,
which implies that $a + \phi(b) = 0$,
whence $a = 0$ and $\phi(b) = 0$.
Since $\phi$ is quadratic, the set 
$\{ b \in  B \mid \xi \in \stab(b) \ \mbox{and}  \phi(b) = 0 \}$
is connected.

\end{proof}
\section{Properties of the energy functionals}

In this section, we will put together the results of the previous  sections to
prove that the energy functions $f_n \colon Y_n \to \R$ are
self-perfecting Morse-Kirwan functions.
The main results of this chapter are summarized in the
following two propositions.

\begin{Proposition}\labell{fnisMK}
The functions $f_n \colon Y_n \to \R$ are Morse in the sense of Kirwan. 
\end{Proposition}

\begin{Proposition}\labell{fnfixed}
Let $C$ be a component of the critical set of $f_n$
and let $E_C^-$ be the negative normal bundle at $C$.
Then there exists a subtorus $T \subset G$ and a $Z(T)$ invariant
subset $B \subset C^T$
so that the natural map $G \times_{Z(T)} B \to C$ is an
equivariant homeomorphism.
Moreover,
$(E_C^-)^T$ is a subset of the zero section of $E_C^-$.
Here, $Z(T)$ denotes the centralizer of $T$.
\end{Proposition}

We now concentrate on proving these Propositions.

\begin{Lemma}\labell{critsetAMM}
The function $f_n \colon Y_n \to \R$ is critical at a point
$y = (g_1,\dots,g_{n}, m) \in Y_n$ if
and if only if there exists $\xi \in \fg$ 
with $|\xi| < n \rho$ such that: 
\begin{enumerate}
\item $g_i = \exp(\frac{i}{n} \xi)$ for all $i$ and
\item   $m \in M^\xi$.
\end{enumerate}
\end{Lemma}

\begin{proof}
This follows immediately from
Lemma \ref{critset},  Lemma \ref{firstorder}, and   the
fact that a path $\gamma \colon [0,1] \to G$ with $\gamma(0) = e$
is a geodesic  if and only
there exists $\xi \in \fg$ so that $\gamma(t) = \exp(t \xi)$.
\end{proof}

Lemma \ref{critsetAMM} and Proposition \ref{goodcriticalsets}
give the following explicit description
of the critical points of $f_n$.

\begin{corollary}\labell{finitecritical}
The critical set of the function $f_n$ has a finite number
of components.
Each connected component $C$ of the critical set of $f_n$
is the $G$-orbit of a connected component of the set
\begin{equation}\labell{crit}
C_\xi := \left. \left\{ \left(
\exp\left( \frac{\xi}{n} \right),\ldots,
\exp\left( \frac{(n-1)  \xi}{n}\right),\exp(\xi),m \right)
\subset Y_n \ \right| \ m \in M^\xi  \right\}
\end{equation}
for some  $\xi \in \fg$.
\end{corollary}

Note that the components of the critical set of $f_n$ need not be manifolds.

\begin{Definition}\labell{minman}
Let $C$ be a connected component of the set $G\cdot C_\xi$ for some $\xi$.
For a sufficiently small neighborhood $U\subset Y_n$ of $C$,
define $\Sigma_C = (G \cdot Y_n^\xi)  \cap U$.
\end{Definition}

\begin{Lemma} \labell{Sigma is manifold}
Let $C$ be a component of the critical set of $f_n$.  Then  $\Sigma_C$
is a locally closed submanifold of $Y_n$.

Moreover, there exists $\xi \in  G$
and a connected component
$B$ of $C_\xi$ so that the natural map
$G \times_{\Stab(\xi)} B \to C$ is
an equivariant homeomorphism.
Let $N(\Sigma_C)$ denote
the normal bundle to $\Sigma_C.$
Then $(N(\Sigma_C))^\xi$ is a subset
of the zero section of $N(\Sigma_C)$.
\end{Lemma}

\begin{proof}
Let $c \in C$; then $c \in C_\xi$ for some $\xi$.
Let $B$ the connected component of $C_\xi$ containing
$c.$
First, note that $Y_n^\xi$ is a manifold, since it is
the fixed point set of a subgroup of $G$.
The stabilizer of $\exp(\frac{\xi}{n})$ is $\Stab(\xi)$.
(Recall that $Y_n$ is constructed so that 
$\exp(\frac{\xi}{n})$ is  not conjugate to $e$ along the geodesic
$\exp(t\xi)$.)
Thus for every $x \in B$, and hence all nearby points in $Y_n^\xi$,
the stabilizer of $x$ is contained
in $\Stab(\xi)$.
Moreover, the action of $\Stab(\xi)$ takes $Y_n^\xi$ to itself.
Therefore, $G \cdot Y_n^\xi = G \times_{\Stab(\xi)} Y_n^\xi$
is a manifold near $C$.
Similarly, by Corollary \ref{finitecritical},
$C$ is homeomorphic to $G \times_{\Stab(\xi)} B$.
Finally, it is clear that $\exp( t \xi)$ acts
freely on the complement of the zero section
of the normal bundle to $Y^\xi_n$, and hence
also on the complement of the zero section of 
the normal bundle to $\Sigma_C$.
\end{proof}

To finish our proof, we need the following fact about Lie groups.

\begin{Lemma}\labell{factaboutjacobi}
Let $G$ be a compact Lie group.
Let $W$ be a Jacobi field along
the geodesic $t \rightarrow \exp(\xi t)$ for some $\xi \in \fg$.
If $W(0) = 0$ and $[\frac{DW}{Dt}(0),\xi] = 0$,
then  $[W(t),\xi] = 0$ for all $t$.
Moreover, $W(1)$ is a non-zero
multiple of $\frac{DW}{Dt}(1)$.
\end{Lemma}

\begin{proof}
Given vector fields $X$, $Y$, and $Z$ on $G$, let $R(X,Y) Z$
denote the curvature tensor.
Fix $\xi \in \fg$.
Let $\gamma$ be the geodesic $\gamma(t) = \exp(t \xi)$.
Define a linear transformation $K \colon T_e G \to T_e G$ by
$K(W) = R(\xi,W) \xi$.
Since $K$ is self-adjoint, we may choose an orthonormal basis $U_1, \ldots
U_n$ for $T_eG$ so that $K(U_i) = e_i U_i$, where $e_i$ are the
eigenvalues of $K$.
We can now extend the $U_i$ to vector fields along $\gamma$ by
parallel translation.
Any Jacobi field along $\gamma$ which vanishes at $0$ can now
be written uniquely as
$W = \sum c_i w_i U_i$, where
$w_i = t$ if $e_i = 0$, $w_i = \sin ( \sqrt{e_i} t)$ if $e_i  > 0$,
and $w_i = \sinh (\sqrt{-e_i} t)$ if $e_i < 0$
(See \cite{Milnor}).
From this, it is clear that $DW/Dt(0) = \sum c_i d_i U_i$,
where $d_i \neq 0$ for all $i$.

Now, since $R(X,Y)V = 1/4 [[X,Y],V]$ for left invariant
vector field, it follows that if $[U_i,\xi]= 0$, then
$e_i = 0$.
So we have $W(t) = \sum_i c_i t U_i$ and $DW/Dt(0) = \sum c_i U_i$.
So, in fact $W(1) = DW/Dt(0)$.
\end{proof}

\begin{Lemma} \labell{C is min}
Let $C$ be any connected component of the critical set
of $f_n$.
Then $C$ is the subset of $\Sigma_C$
on which $f_n$ takes its minimum value.
\end{Lemma}

\begin{proof}
Recall that $\Sigma_C = G \cdot (Y^\xi_n \cap U)$, for a sufficiently
small neighborhood $U$ of $C_\xi$, for some $\xi \in \fg$.

Suppose that $f_n$ takes its minimal value at $x = (g_1,\ldots,g_n,m) \in
Y^\xi_n \cap U \subset (\Stab(\xi))^n \times M^\xi$.
We will show $x \in C_\xi$.

Since $x$ is near $C_\xi$ and $\Phi(m) \in \Stab(\xi)$,
there exists $\eta \in \stab(\xi)$ near $\xi$ such that
$\exp(\eta) = \Phi(m)$.  Since $\xi$ is in the center of $\Stab(\xi)$,
and $\eta$ is nearby, there are no conjugate points to $e$ in $\Stab(\xi)$
along the geodesic $t \to \exp(t\eta)$.  So the shortest path from $e$ to
$\exp(\eta)$ near $t \to \exp(t\xi)$ is the geodesic $t \to \exp(t \eta)$.
Hence, we may assume $x = (\exp(\eta/n), \exp(2\eta/n),\ldots,
\exp(\eta),m)$.

We now compute $f_n(x)$ using
the normal form theorem.
Let $K=\Stab(\zeta)$, and let $H=\Stab(m)$.
By Proposition \ref{qnormalform}, a neighbourhood
of $m$ in $M$ is given by a neighbourhood of $(e,0,0)$ in
$Y=G\times_K((\fh^\perp\cap \fk) \times V)$; under this identification
the moment map $\Phi$ is identified with the map
$\Phi \colon Y \to G$  given by
$$\Phi([g,\gamma,v]) = Ad_g(\zeta\exp(\gamma + \phi(v))),$$
\noindent where $\phi \colon V \to \fh^*\simeq \fh$
is the moment map for the linear $H$ action
on the symplectic vector space $V$.

\mute{this line above used to say ``K action'' but I think we mean H.}

In terms of these coordinates, $f_n(x)$ is given by the
the square of the length of the geodesic $\exp (t\eta)$ from
$e$ to $\exp(\eta)=\zeta\exp((\gamma + \phi(v))$, that is, by
$|\eta|^2$.
But
$\eta \in Lie(C(\xi))$, so that $\zeta^{-1}\exp(\eta) = \exp (\eta-\xi)$.
Thus $f_n(x)=|\eta|^2 = |\xi+ \gamma + \phi(v)|^2$.
Since the image of $((\fh^\perp\cap\fk) \times V^\xi)$ under
$\gamma + \phi(v)$ is a hyperplane in $\fg$
through the origin and perpendicular to $\xi$, the distance from
that hyperplane to $-\xi$ is minimized at the origin.
It follows
that the function $f_n(x)$ attains its minimum where $\gamma+\phi(v)=0$,
that is, where $\Phi(m)=\zeta$.
\end{proof}

\begin{Lemma} \labell{Hessian is nondeg}
Let $C$ be any connected component of the  critical set.
For every point $x \in C$, the
null-space of the Hessian $Hf_n $ is contained in the
tangent space to  $\Sigma_C$.
\end{Lemma}

\begin{proof}
Let  $y = (g_1,\ldots,g_n,m) \in Y_n$ be a critical point of $f_n$.
Then there exists $\xi\in\fg$ such that $g_i= \exp(i\xi/n)$ and
$m\in M^\xi$.  Denote by $\lambda$ the geodesic $t\to \exp(t\xi)$.
%Recall that we have denoted by $\lambda_y$ the
%path $t\to \exp(t\xi)$ associated to $y$.
Let  $\eta = (X_1,\ldots,X_n,Y)$ be in the null space of the Hessian $H f_n$.
We need to show that $\eta$ is
in the tangent bundle to $\Sigma_C$.

Let  $W$ be the unique
Jacobi field on $\lambda$ which is smooth except (possibly) at
$g_i$, and such that $W(0) =  0$ and $W(i\xi/n) = X_i$
for all $1 < i \leq n$.
Let $\eta' = (\Delta_1 \frac{DW}{Dt},\ldots,\Delta_{n-1} \frac{DW}{Dt},0,0)$.
Then
$H f_n (\eta, \eta') = \sum_{i < n} | \Delta_i \frac{D W}{D t}|^2 $.
Hence, $W$ must be smooth.

Let  $\alpha = DW/Dt(0) \in \fg$.
The exact sequence
$0\to {\rm ker~}([\xi,\cdot]) \to \fg \to [\xi,\fg]\to 0$
shows that $\fg = {\rm ker~}([\xi,\cdot]) \oplus [\xi,\fg]$.
Hence $\alpha$ can be written as a sum
$\alpha = \beta + [\gamma,\xi]$, where $[\beta,\xi] = 0$.

Let $W_\gamma$ be the unique
unbroken Jacobi field so that $W_\gamma(0) = 0$ and
$\frac{DW_\gamma}{Dt}(0) = [\gamma,\xi]$.
Then $W_\gamma$ is the Jacobi field on the geodesic path
$t\to \exp(t\xi)$ induced by conjugation by $\exp(s\gamma)$.
Let $\hat{\gamma} \in T_m M$ be the value at $m$
of the vector field $\gamma^M$ on $M$ induced by $\gamma$.
Then, by Lemma \ref{invar},
$$\left<DW_\gamma/Dt(1),\Phi_*(Y')\right> + H^\xi \Phi(\hat{\gamma},Y')=0$$ 
\noindent for all $Y' \in T_m M$.
Let $\eta_\gamma = (W_\gamma(1/n),\cdots,W_\gamma(1),\hat{\gamma})$.
Then $\eta_\gamma$
is in the null-space of the Hessian.
Therefore, $ \eta - \eta_\gamma$ is also in the null-space of the Hessian.
Moreover,  since $\eta_\gamma$  is the value at $y$ of the vector field on
$G^m \times M$ induced by $\gamma$, $\eta_\gamma$ is in the tangent space of $\Sigma_C$.
Thus, to show that $\eta$ is
an element of $T_y \Sigma_C,$ it suffices to show that $\eta - \eta_\gamma$ is in the tangent
space of $\Sigma_C$.
Equivalently, we can assume that $\eta$ is such
that $[DW/Dt(0),\xi] = 0$.

But by Lemma \ref{factaboutjacobi},
if $DW/Dt(0)$ commutes with $\xi$, so does $W(1)$, and
moreover, $ c DW/Dt(1) = W(1)$ for some non-zero real number $c$.
By assumption,  $\Phi_*(Y) = W(1)$.
So by Lemma \ref{hamimage}
there also exists $Y' \in T_m M^\xi$ so that $\Phi_*(Y') = W(1)$.
But then $H^\xi \Phi(Y,Y') = 0$ because by Proposition \ref{commuting}
$Y'$ is in the null-space of
$H^\xi \Phi$.  Therefore, letting $\eta'=(X_1,\cdots,X_n,Y')$, we have

 $$Hf_n(\eta,\eta')= \left< DW/Dt(1), \Phi_*(Y') \right> =
 \left< DW/Dt(1), W(1) \right>  = c |DW/Dt(1)|^2. $$
But our assumption is that $\eta$ is in the null space of the
Hessian $H f_n$.  So we have $DW/Dt(1) = 0$ and $\eta=(0,\cdots,0,Y)$
where $Y \in {\rm ker~}{(\Phi_*)}_y.$

This means that for any $\eta'=(X_1',\cdots,X_n',Y')$,
$Hf_n(\eta,\eta') =  H^\xi\Phi(Y,Y')=0$.
Thus $Y$ is in the nullspace of the Hessian $H^\xi\Phi_y$.
By Proposition \ref{nondeg}, $Y\in T_m(M^\xi),$
so that $\eta\in T_y\Sigma_C$, as needed.

\end{proof}

\begin{proof}[Proof of Propositions \ref{fnisMK} and \ref{fnfixed}]

The set of critical points
of $f_n$ has a finite number of components 
by Corollary \ref{finitecritical}. 
By Lemma \ref{Sigma is manifold},
$\Sigma_C$ is a locally closed
submanifold of $Y_n$.
Moreover, there exists $\xi \in \fg$ and a connected component
$B$ of $C_\xi$ so that the natural map
$G \times_{\Stab(\xi)} B \to C$
is an equivariant homeomorphism.
Let $T$ be the closure of the 
one-parameter subgroup $\exp(\xi t)$ generated by $\xi$.
By Lemma \ref{Sigma is manifold},
$(N(\Sigma_C))^T$ is a subset of 
the zero section of the bundle
$N(\Sigma_C)$. 
The $T$-action
also endows $N(\Sigma_C)|_{B}$
with an orientation; 
this orientation extends to all of $N(\Sigma_C)$.
Finally the submanifold $\Sigma_C$ satisfies
conditions 2(a) and 2(b) of Definition \ref{morsekirwan}
by Lemma \ref{C is min} and Lemma \ref{Hessian is nondeg}.
This proves Proposition \ref{fnisMK}.
Since $E_C^-$ is a subbundle 
of $N(\Sigma_C)|_C$, Lemma \ref{Sigma is manifold}
also  proves Proposition \ref{fnfixed}.
\end{proof}

\section{Kirwan's extension of Morse theory}\labell{kirsec}

In this section we give a brief description of Kirwan's extension
of Morse theory as given in \cite{kirwan}.

\begin{Definition}\labell{morsekirwan}
A smooth function $f \colon X \to \R$ is
{\bf Morse in the sense of Kirwan} if
\begin{enumerate}
\item
The set of critical points of $f$ has a finite number of
components.
\item
For each component $C$ of the critical set
there exists a locally closed submanifold $\Sigma_C$
with an orientable normal bundle in $X$ such that:
\begin{enumerate}
\item $C$ is is the subset of $\Sigma_C $
on which $f$ takes its minimum
value.
\item At every point $x \in C$, every element $\eta \in T_x M$ that
lies in the null space of the Hessian lies in
 the tangent space $T_x \Sigma_C$.
\end{enumerate}
Given a component $C$ of the critical set, we
call $\Sigma_C$ a {\bf minimizing manifold} for $f$ along $C$.
\end{enumerate}\end{Definition}

\begin{Definition}\labell{negnor}
Let $f$ be Morse in the sense of Kirwan.  Given any critical
set $C$ of $f$, and any point $p \in C$, the Hessian $Hf_p$
splits the normal bundle to $\nu(\Sigma_C)$ into positive,
negative, and null eigenspaces.  We denote by $E_C^-$, the {\bf negative
normal bundle at $C$}, the
vector bundle on $C$ given by the negative eigenspaces of $Hf$,
and define the {\bf index} $\lambda_C$ of $C$ as the
dimension of $E_C^-$.
\end{Definition}

This definition of a Morse function in the sense of Kirwan
is not identical
to the definition given in \cite{kirwan}.
According to Kirwan, a function
$f \colon X \to \R$ is {\bf minimally degenerate}
if
\begin{enumerate}
\item
The set of critical points of $f$ is a finite union
of disjoint closed subsets $C$ of $X$ on each of which $f$ takes a constant
value.
\item
For each component $C$ of the critical set
there exists a locally closed submanifold $\Upsilon_C$
with an orientable normal bundle in $X$ such that:
\begin{enumerate}
\item $C$ is the subset of $\Upsilon_C$ on which $f$ takes its minimum
value.
\item At every point $x \in C$, the tangent space $T_x\Upsilon_C$
is maximal among those subspaces of $T_x X$ where the Hessian
$Hf_x$ is positive semidefinite.
\end{enumerate}
\end{enumerate}

However, we have the following Lemma:

\begin{Lemma}
If a smooth function $f \colon X\to \R$ is Morse in the sense of Kirwan,
it is minimally degenerate.
\end{Lemma}

\begin{proof}  
Suppose $f$ is Morse in the sense of Kirwan.  Let $C$ be
a critical set of $f$ and let $\Sigma_C$ be the corresponding
minimizing manifold.  Let $U_C$ be a tubular neighbourhood
of $\Sigma_C$ and identify $U_C$ with $N\Sigma_C$.  Let $\Upsilon_C$
be the positive normal bundle to $\Sigma_C$, considered
as a submanifold of $U_C$.  It is clear
that for each $x\in C$, $T_x\Upsilon_C$ is maximal among subspaces
where $Hf_x$ is positive semidefinite.
\end{proof}

This lemma allows us to apply the results of Chapter 10 of
\cite{kirwan}, which constructs a version of Morse theory
for minimally degenerate functions.  We may summarize the
results of this construction in the following

\begin{Proposition}\labell{morsekirwanlemma}
Let $f$ be a minimally degenerate function on a smooth
manifold $X$, and let $c$ be a critical value of $f$,
corresponding to a critical set $C$.  For $\epsilon$
sufficiently small, define  $$M_\pm = f^{-1}(-\infty, c\pm \epsilon).$$

Denote by $D^-$ the negative disc bundle to $\Sigma_C$, restricted
to $C$, and by $S^-$ the corresponding sphere bundle.
Then there exists a long exact sequence

\begin{equation}\labell{exactseq}
\vcenter{\xymatrix{
\cdots
\ar[r] &
H^*(D^-,S^-)
\ar[r]\ar[d]^\simeq &
H^*(M_+)
\ar[r] \ar[d] &
H^*(M_-)\ar[r] &  \cdots
\\
&
H^{*-{\rm dim} D^-} (C)
\ar[r]^{~~\cup e(D^-)}  &
H^*(C) &
\\
}}
\end{equation}

\noindent where $e(D^-)$ is the Euler class of the bundle $D^-$.
If a compact group $G$ acts on $X$, and the function $f$ is invariant,
the same result holds in equivariant cohomology.
\end{Proposition}


\begin{thebibliography}{DGMW}

\bibitem[AMM]{AMM} A. Alexeev, A. Malkin, and E. Meinrenken, Lie group
valued moment maps.  {\em J. Diff. Geom.} {\bf 48} 445-495 (1998)

\bibitem[AB]{ab} M. Atiyah, R. Bott.  The Yang-Mills
functional over a Riemann surface.  {\em Phil. Trans. Roy.
Soc.} {\bf A308}, 523-615 (1982)

\bibitem[B]{Bott} R. Bott, The stable homotopy of the
classical groups. {\em Ann. Math.} {\bf 70}, 313-337 (1957)

\bibitem[BT]{bt} R. Bott, L. Tu.  {\em Differential forms
in algebraic topology}.  Springer Verlag, 1982.

\bibitem[CM]{cm} A.L. Carey and M. K. Murray.  String structures
and the path fibration of a group.  {\em Commun. Math. Phys.}
{\bf 141}, 441-452 (1991)

\bibitem[F]{frankel} T. Frankel, Fixed points on Kahler
manifolds.  {\em Ann. Math.} {\bf 70}, 1-8 (1959)

\bibitem[GS]{gs} V. Guillemin and S. Sternberg.  {\em Symplectic
Techniques in Physics.}  Cambridge University Press, 1984.

\bibitem[GS2]{gs2} V. Guillemin and S. Sternberg.  {\em Supersymmetry
and equivariant de Rham theory.}  Springer Verlag, 1999.

\bibitem[Kil]{killingback}  T. Killingback.  World-sheet anomalies
and loop geometry.  {\em Nucl. Phys.} {\bf B288}, 578-588 (1987)

\bibitem[Kir]{kirwan} F. C. Kirwan,  {\em The cohomology of quotients
in symplectic and algebraic geometry}, Princeton University Press, 1984.

\bibitem[M]{Milnor} J. Milnor, {\em Morse Theory}.  Princeton University
Press, 1963

\bibitem[PS]{ps} A. Pressley, G. Segal.  {\em Loop groups}.
Oxford University Press, 1986.

\bibitem[W]{W} Chris Woodward, {\em A Kirwan-Ness stratification for
loop groups}

\end{thebibliography}
\end{document}